\documentclass[a4paper,dvipsnames,fleqn]{article}

\usepackage{amsmath}
\usepackage{amssymb}
\usepackage{amsthm}
\usepackage{bm}
\usepackage[backend=biber,bibencoding=utf8,bibstyle=alphabetic,citestyle=alphabetic,sorting=nyt,maxnames=4]{biblatex}
	\AtEveryBibitem{\clearfield{doi}}
	\AtEveryBibitem{\clearfield{isbn}}
	\AtEveryBibitem{\clearfield{issn}}
	\AtEveryBibitem{\clearfield{url}}
	
	\renewbibmacro{in:}{\ifentrytype{article}{}{\printtext{\bibstring{in}\intitlepunct}}}
	\AtEveryBibitem{\ifentrytype{book}{\clearfield{pages}}{}}
\usepackage[colorlinks=true,allcolors=purple]{hyperref}
\usepackage{cleveref}
	\crefformat{equation}{#2(#1)#3}
\usepackage{dsfont}
\usepackage{enumitem}
\usepackage[OT2,OT1]{fontenc}
\usepackage[cal=cm,scr=boondoxo]{mathalfa}
\usepackage{mathtools}
\usepackage{pifont}
\usepackage{stmaryrd}
\usepackage{textcomp}
\usepackage{tikz}
	\usetikzlibrary{arrows}
	\usetikzlibrary{patterns}
\usepackage{tikz-cd}
\usepackage[textwidth=3cm,colorinlistoftodos]{todonotes}
\usepackage{xfrac}

\numberwithin{equation}{section}			% Enumeration of equations
\swapnumbers						% Enumeration of theorem-like environments

\newcommand\cyr{%					% cyrillic font
\renewcommand\rmdefault{wncyr}%
\renewcommand\sfdefault{wncyss}%
\renewcommand\encodingdefault{OT2}%
\normalfont
\selectfont}
\DeclareTextFontCommand{\textcyr}{\cyr}

\makeatletter						% content of ctan-package "capt-of"
\newcommand\captionof[1]{\def\@captype{#1}\caption}	
\makeatother

\newenvironment{Enumerate}{\begin{enumerate}[label={\rm({\roman*})}]}{\end{enumerate}}

\newcommand{\descriptionlabelsave}{}		% Itemized list
\newenvironment{Itemize}{%
	\renewcommand{\descriptionlabelsave}{\descriptionlabel}\renewcommand{\descriptionlabel}{$\triangleright$}%
	\begin{description}[leftmargin=15pt,itemindent=-5.2pt]}{%
	\end{description}\renewcommand{\descriptionlabel}{\descriptionlabelsave}}

\newcounter{StepsCount}				% Enumerated list with no indentation (e.g. steps in proof)
\newenvironment{Elist}{%
	\begin{list}{\ding{\value{StepsCount}}}{\usecounter{StepsCount} \leftmargin=0pt \labelwidth=12pt \itemindent=\labelwidth%
	\itemsep=5pt\listparindent=\parindent} \setcounter{StepsCount}{191}}{\end{list}}
\newcounter{StepsRefCount}

\newenvironment{Ilist}{%			% Itemised list without indentation (not-enumerated steps)
	\begin{list}{$\triangleright$}{\leftmargin=0pt \labelwidth=11pt \itemindent=\labelwidth%
	\itemsep=5pt\listparindent=\parindent}}{\end{list}}

\theoremstyle{plain}
	\newtheorem{lemma}{Lemma}[section]
	\newtheorem{proposition}[lemma]{Proposition}
	\newtheorem{theorem}[lemma]{Theorem}
	\newtheorem{corollary}[lemma]{Corollary}
	\newcommand{\GenericTheoremName}{}\newtheorem{generictheorem}[lemma]{\GenericTheoremName}
\theoremstyle{definition}
	\newtheorem{definition}[lemma]{Definition}
	\newcommand{\GenericDefinitionName}{}\newtheorem{genericdefinition}[lemma]{\GenericDefinitionName}
\theoremstyle{remark}
	\newtheorem{remark}[lemma]{Remark}
	\newtheorem{example}[lemma]{Example}
	\newcommand{\GenericRemarkName}{}\newtheorem{genericremark}[lemma]{\GenericRemarkName}
\newenvironment{Theorem}{\begin{theorem}}{\par\noindent\rule{5em}{1pt}\end{theorem}}%6.4
\newenvironment{Corollary}{\begin{corollary}}{\par\noindent\rule{5em}{1pt}\end{corollary}}%6.5
\newenvironment{Definition}{\begin{definition}}{\par\noindent\rule{5em}{1pt}\end{definition}}%7
\newenvironment{Remark}{\begin{remark}}{\par\noindent\rule{5em}{0.5pt}\end{remark}}%5.2
\newcommand{\mc}[1]{{\mathcal{#1}}}			% --- abbreviation ---
\newcommand{\ms}[1]{{\mathscr{#1}}}			% --- abbreviation ---
\newcommand{\bb}[1]{{\mathbb{#1}}}			% --- abbreviation ---
\newcommand{\ov}{\overline}				% --- abbreviation ---
\newcommand{\mr}{\mathring}				% --- abbreviation ---
\DeclareMathOperator{\RE}{Re}				% real part
\renewcommand{\Re}{\RE}
\DeclareMathOperator{\IM}{Im}				% imaginary part
\renewcommand{\Im}{\IM}

\newcommand{\Side}[1]{\hfill{#1}\kern10pt}		% text put on the left side of line with offset

\newcommand{\smmatrix}[4]{\Bigl(			% small matrix for use in textline
\begin{smallmatrix}
\hspace*{-0.2ex} #1 \hspace*{0.2ex} & \hspace*{0.2ex} #2 \hspace*{-0.2ex}
\\[0.5ex]
\hspace*{-0.2ex} #3 \hspace*{0.2ex} & \hspace*{0.2ex} #4 \hspace*{-0.2ex}
\end{smallmatrix}
\Bigr)}
\newcommand{\smfrac}[2]{{\textstyle\frac{{#1}}{{#2}}}}	% small fraction for use in displayed formulas

\newcommand{\Dummy}{\text{\textvisiblespace\kern1pt}}	% Platzhaltersymbol fuer Funktionsargumente
\newcommand{\Smallo}{{\rm o}}				% small o
\newcommand{\BigO}{{\rm O}}				% big o

\DeclareMathOperator{\Id}{id}				% identity map
\DeclareMathOperator{\Span}{span}			% linear span
\DeclareMathOperator{\Cls}{cls}				% closed linear span
\DeclareMathOperator{\Tr}{tr}				% trace
\DeclareMathOperator{\SL}{SL}				% special linear group
\DeclareMathOperator{\Clos}{Clos}			% closure
\DeclareMathOperator{\Rank}{rank}			% rank of a matrix
\DeclareMathOperator{\GL}{GL}				% general linear group

\newcommand{\DS}{\mid\mkern3mu}				% delimiter for set definition
\newcommand{\DQ}{\mkern6mu}				% distance for successive quantors
\newcommand{\DP}{{\colon\kern5pt}}			% delimiter for predicate formula
\newcommand{\DF}{\colon}				% delimiter for function domain/codomain
\newcommand{\DE}{\mathrel{\mathop:}=}			% defining equality
\newcommand{\ED}{=\mathrel{\mathop:}}			% defining equality
\newcommand{\DI}{\mathrel{\mathop:}\Leftrightarrow}	% defining equivalence
\newcommand{\DD}{\mkern4mu\mathrm{d}}			% distance and rm-d for integration differential
\newcommand{\CAS}{&\text{if}\ }				% delimiter in "cases" environment
\newcommand{\CASO}{&\text{otherwise}}			% delimiter in "cases" environment

\newcommand{\HB}{\mc H\mkern-2mu\mc B}			% Hermite-Biehler class
\DeclareMathOperator{\Hol}{Hol}				% Holomorphic functions
\DeclareMathOperator{\Sub}{Sub}				% Chain of subspaces
\newcommand{\PW}{{\mc P\hspace*{-1pt}W\!}}		% Paley-Wiener spaces

\begin{document}

\begin{flushleft}
	{\Large\bf Homogeneous spaces of entire functions}
	\\[5mm]
	\textsc{
	Benjamin Eichinger
	\,\ $\ast$\,\ 
	Harald Woracek
		\hspace*{-14pt}
		\renewcommand{\thefootnote}{\fnsymbol{footnote}}
		\setcounter{footnote}{2}
		\footnote{The first author was supported by the stand alone project P~33885 of the Austrian Science Fund (FWF).
			The second author was supported by the joint project I 4600 of the Austrian Science Found 
			(FWF) and the Russian foundation of basic research (RFBR).}
		\renewcommand{\thefootnote}{\arabic{footnote}}
		\setcounter{footnote}{0}
	}
	\\[6mm]
	{\small
	\textbf{Abstract:}
	Homogeneous spaces are de~Branges' Hilbert spaces of entire functions with the property that certain 
	weighted rescaling transforms induce isometries of the space into itself. A classical example of a homogeneous space is
	the Paley-Wiener space of entire functions with exponential type at most $a$ being square integrable on the real axis.
	Other examples occur in the theory of the Bessel equation. Being homogeneous is a strong property,
	and one can describe all homogeneous spaces, their structure Hamiltonians, and the measures associated with chains of 
	such spaces, explicitly in terms of powers, logarithms, and confluent hypergeometric functions. 
	\\[2mm]
	The theory of homogeneous spaces was in large parts settled by L.de~Branges in the early 1960's. 
	However, in his work some 
	connections and explicit formulae are not given, some results are stated without a proof, and last but not least a 
	mistake occurs which seemingly remained unnoticed up to the day. 
	\\[2mm]
	In this paper we give a detailed account on homogeneous spaces. We provide explicit proofs for all formulae
	and relations between the mentioned objects, and correct the mentioned mistake. 
	\\[3mm]
	}
\end{flushleft}

%%%%%%%%%%%%%%%%%%%%%%%%%%%%%%%%%%%%%%%%%%%%%%%%%%%%%%%%%%%%%%%%%%%%%%%%%%%%%%%%%%%%%%%%%
% disable in final version
%\pagenumbering{roman}
%\fbox{
%\parbox{100mm}{
%\hspace*{0pt}\\[1mm]
%\centerline{{\Large\ding{45}}\quad\,{\large\sc Draft}\quad{\Large\ding{45}}}
%\hspace*{0pt}\\[-2mm]
%\textcircledP\ \ Preliminary version Sat 26 Sep 2026 10:23
%\\[2mm]
%\ding{233}\quad Use pdflatex/biber to compile
%\\[-1mm]
%}
%}
%\tableofcontents
%\listoftodos
%\newpage
%\pagenumbering{arabic}
%\setcounter{page}{1}
%
%%%%%%%%%%%%%%%%%%%%%%%%%%%%%%%%%%%%%%%%%%%%%%%%%%%%%%%%%%%%%%%%%%%%%%%%%%%%%%%%%%%%%%%%%

%---------
%   TEXTBODY
%---------

%
%
%
\section[\textcolor{ForestGreen}{Introduction}]{Introduction}
%ForestGreen / Dandelion / BrickRed / Fuchsia
%
%

The theory of \emph{Hilbert spaces of entire functions} was founded by L.de~Branges in the late 1950's and further developed in a
series of papers. A comprehensive presentation appeared as the monograph \cite{debranges:1968}. These spaces are reproducing
kernel Hilbert spaces whose elements are entire functions and which possess certain additional properties.

\begin{Definition}
\label{M101}
	A \emph{de~Branges space} is a Hilbert space $\mc H$ which satisfies 
	the following axioms.
	\begin{Itemize}
	\item The elements of $\mc H$ are entire functions and
		for each $w\in\bb C\setminus\bb R$ the point evaluation functional 
		$F\mapsto F(w)$, $F\in\mc H$, is linear and continuous in the norm $\|.\|_{\mc H}$ of $\mc H$.
	\item For each $F\in\mc H$, also the function $F^\#(z)\DE\ov{F(\ov z)}$
		belongs to $\mc H$ and $\|F^\#\|_{\mc H}=\|F\|_{\mc H}$.
	\item If $w\in\bb C\setminus\bb R$ and $F\in\mc H$ with $F(w)=0$, then
		\[
			\frac{F(z)}{z-w}\in\mc H
			\quad\text{and}\quad
			\Big\|\frac{z-\ov w}{z-w}F(z)\Big\|_{\mc H}=\big\|F\big\|_{\mc H}
			.
		\]
	\end{Itemize}
	Throughout this paper the notion of a de~Branges space shall additionally
	include the following requirement:
	\begin{Itemize}
	\item For each $t\in\bb R$ there exists $F\in\mc H$ with $F(t)\neq 0$.
	\end{Itemize}
\end{Definition}

\noindent
The first three properties imply that point evaluations are also continuous at each real point $w$, and together with the fourth
property it follows that the third axiom holds also for all real points $w$. 

Hilbert spaces of entire functions are an -- in essence -- equivalent view on entire operators in the sense of M.G.Krein, e.g.\ 
\cite{gorbachuk.gorbachuk:1997}, which makes them relevant in an operator theoretic context. Further, they are equivalent to
certain shift-coinvariant subspaces of the Hardy space, which makes them relevant in a function theoretic context. 
Typical areas where de~Branges spaces occur are the spectral theory of Sturm-Liouville or Krein-Feller operators, 
e.g.\ \cite{remling:2002,kaltenbaeck.winkler.woracek:bimmel}, 
models for symmetric operators \cite{martin:2011,aleman.martin.ross:2013}, 
interpolation and sampling e.g.\ \cite{ortega.seip:2002,baranov.belov:2011}, or Beurling-Malliavin type theorems 
e.g.\ \cite{havin.mashreghi:2003a,havin.mashreghi:2003b}. 
Our standard reference is \cite{debranges:1968}; other references are e.g.\ 
\cite{baranov.woracek:SprRef-DB,romanov:1408.6022v1,remling:2018}.

One can think of de~Branges' theory as a generalisation of Fourier analysis. The maybe best known examples 
of de~Branges spaces are the Paley-Wiener spaces: given $a>0$ the set 
\[
	\PW_a\DE \big\{F\DS \parbox[t]{90mm}{$F$ entire with exponential type $\leq a\ \wedge$
	$\frac 1\pi\int_{\bb R}|F(x)|^2\DD x<\infty$}\big\}
\]
is a de~Branges space when endowed with the $L^2$-scalar product induced by the Lebesgue measure on $\bb R$.
These spaces, and de~Branges' structure theory applied with them, are closely related to the exponential function and the 
classical sine- and cosine transforms. 

Paley-Wiener spaces have several very specific structural properties. 
One of them is that each space $\PW_a$ is well-behaved with respect to a weighted rescaling transform: 
for all $c\in(0,1]$ the map 
\begin{equation}
\label{M98}
	F(z)\mapsto c^{\frac 12}F(cz)
\end{equation}
maps $\PW_a$ isometrically into itself. 
The same phenomenon occurs in the theory of Bessel functions and, more generally, confluent hypergeometric functions. 
The only difference being that for the spaces related with such functions the power $\frac 12$ in \cref{M98} 
has to be replaced with another power. 
This fact served as a motivation for L.de~Branges to formulate the following axiomatic definition, cf.\ 
\cite[\S50]{debranges:1968}.

\begin{Definition}
\label{M70}
	Let $\nu>-1$ and let $\mc H$ be a de~Branges space. Then $\mc H$ is called \emph{homogeneous of order $\nu$}, 
	if for all $c\in(0,1]$ the map 
	\begin{equation}
	\label{M100}
		F(z)\mapsto c^{\nu+1}F(cz)
	\end{equation}
	maps $\mc H$ isometrically into itself. 
\end{Definition}

\noindent
Being homogeneous is a very strong property, and not many de~Branges spaces possess it. 

The source of homogeneity of a space $\mc H$ can be pinpointed in different ways.
\begin{enumerate}[label={\rm({\arabic*})}]
\item The power series coefficients of a certain entire function determining $\mc H$ satisfy a recurrence relation 
	of a particular form 
	(a linear recurrence for the vector made up from the real- and imaginary parts of the coefficients).
\item The family of reproducing kernels of de~Branges subspaces of $\mc H$ (cf.\ \Cref{M49}) satisfies a functional equation 
	involving a weight and rescaling. 
\item The canonical system given by the chain of de~Branges subspaces of $\mc H$ (\Cref{M25}) has a Hamiltonian of a 
	particular form (involving powers and, possibly, logarithms).
\item The norm of $\mc H$ is the $L^2$-norm given by a measure which is absolutely continuous w.r.t.\ the Lebesgue measure
	and has power density (on the left- and right half-axis separately), and the functions in the space are of bounded type 
	in the upper half-plane.
\end{enumerate}
Most of these topics were investigated in \cite{debranges:1962b} and \cite[Theorem~50]{debranges:1968}, although some of
them are not made explicit. For example it is stated on \cite[p.205]{debranges:1962b} that 
``homogeneous spaces of entire functions are related to Bessel functions and more general confluent hypergeometric functions'', 
but the actual formulae are not given. Also the passage from measures with power density to homogeneous
de~Branges spaces is not treated. 

Homogeneous spaces appeared in some places in the literature. We mention
\cite{holt.vaaler:1996,carneiro.littmann:2014,vaaler:2023} who use the theory in the context of an extremal problem from 
number theory, and 
\cite{goncalves:2017,goncalves.littmann:2018} who study invariance under differentiation and give interpolation formulae. 
It should be noted that in those references only homogeneous spaces occur which enjoy an additional symmetry property and 
are related to Bessel functions. 

In the recent manuscript \cite{eichinger.lukic.woracek:elw-arXiv} homogeneous spaces in their full generality corresponding to
confluent hypergeometric functions play a decisive role.
And it was during writing of that paper, that we found out that \cite{debranges:1962b,debranges:1968}
contains a mistake. In fact, in the
description of all homogeneous spaces of order $-\frac 12$ (which contains the Paley-Wiener case), a 
whole $1$-parameter family of spaces was forgotten; details are explained in \Cref{M17}. 
It seems that previously this mistake was not noticed; fortunately it also does not affect the earlier 
literature mentioned above due to the symmetry present in the spaces of those papers. 

The purpose of our present paper is twofold.
One, we provide a detailed account on homogeneity in de~Branges spaces and all the viewpoints listed above, 
round off the picture by taking some slightly more general or more systematic viewpoints, 
and discuss side results which were not touched upon in de~Branges' original work. 
Two, we correct the mentioned mistake. 
We should say it very clearly that our aim is to provide a comprehensive structured account on homogeneous spaces and 
to make the results accessible also to non-specialists. Hence, fully elaborated proofs of all assertions are included, 
also when some of them simply follow what was done by de~Branges.

Let us briefly describe the structuring of the paper. 
In the second part of this introduction we set up our notation concerning de~Branges spaces and recall several facts which are
needed in order to make the presentation self-contained.
Then we introduce the main players of the paper and provide some of their elementary properties (Section~2). 
There follow two sections (Sections~3 and 4) where we investigate the solutions of the canonical systems occurring 
in the present context: first, approaching them via the power series coefficient recurrence (\Cref{M53})
and, second, giving explicit formulae in terms of special functions (\Cref{M36}). 
In Section~5 we systematically investigate the group action of power-weighted rescalings in the context of de~Branges spaces. 
Putting together all those results, this culminates in Sections~6 and 7 in a complete description of homogeneous spaces.
We determine the structure of their chain of subspaces (\Cref{M16}), and the measures associated with 
chains of homogeneous spaces (\Cref{M37}). 

\subsection[\textcolor{ForestGreen}{De~Branges spaces, chains, and measures}]{De~Branges spaces, chains, and measures}
%ForestGreen / Dandelion / BrickRed / Fuchsia
%

We recall the basics of the theory of Hilbert spaces of entire functions. This compilation is extracted almost exclusively from 
\cite{debranges:1968}. 

\begin{center}
	\textnumero 1
	\\
	\emph{De~Branges spaces via Hermite-Biehler functions}
\end{center}
In \Cref{M101} we used an axiomatic way to introduce de~Branges spaces. On a more concrete level, these objects may also be
introduced via certain entire functions. The reason being that the reproducing kernel of a de~Branges space has a very 
particular form\footnote{%
	Here we denote by $\bb C_+$ the open upper half-plane, i.e., $\bb C_+\DE\{z\in\bb C\DS\Im z>0\}$.
	}.

\begin{definition}
\label{M47}
	A \emph{Hermite-Biehler function} is an entire function $E$ which satisfies 
	\[
		\forall z\in\bb C_+\DP |E(\ov z)|<|E(z)|
		,
	\]
	and has no real zeroes\footnote{%
		We require absence of real zeroes in order to fit our convention from \Cref{M101} that for all real points a 
		de~Branges space should contain elements which do not vanish at that point. 
		}.

	We denote the set of all Hermite-Biehler functions as $\HB$.
\end{definition}

\noindent
For any entire function $E$ we use the generic notation 
\begin{equation}
\label{M50}
	A\DE\frac 12(E+E^\#),\ B\DE\frac i2(E-E^\#),\qquad E=A-iB
	,
\end{equation}
and denote 
\begin{equation}\label{eq:1}
	K_E(z,w)\DE\frac{B(z)A(\ov w)-B(\ov w)A(z)}{z-\ov w}
	.
\end{equation}
Then an entire function $E$ without real zeroes belongs to $\HB$ if and only if $K_E$ is a positive kernel and is not
identically equal to zero.

The connection between de~Branges spaces and the Hermite-Biehler class can be summarised as follows.

\begin{theorem}
\label{M99}
	\phantom{}
	\begin{Enumerate}
	\item Let $E\in\HB$. Then the reproducing kernel Hilbert space $\mc H(E)$ generated by the positive kernel $K_E$ is a 
		de~Branges space. 
	\item Let $\mc H$ be a de~Branges space. Then there exists $E\in\HB$ such that $\mc H=\mc H(E)$. 

	\item Let $E,\tilde E\in\HB$. Then $\mc H(E)=\mc H(\tilde E)$ if and only if there exists $M\in\SL(2,\bb R)$ such that 
		\[
			(A,B)=(\tilde A,\tilde B)M
			.
		\]
	\end{Enumerate}
\end{theorem}

\noindent
This theorem allows to switch between abstract and concrete levels. While the abstract -- axiomatic -- viewpoint is suitable to 
make the connection with operator theory, the concrete viewpoint allows to invoke classical function theory.

The description of a de~Branges space $\mc H$ via the reproducing kernel $K_E$ is implicit, since it involves a completion
process to pass from the linear span of kernels to the whole space. An explicit description reads as follows.

\begin{theorem}
\label{M20}
	Let $E\in\HB$. Then an entire function $F$ belongs to the space $\mc H(E)$, if and only if 
	\[
		\int\limits_{\bb R}\Big|\frac{F(t)}{E(t)}\Big|^2\DD t<\infty
		\ \wedge\ 
		\forall z\in\bb C\DP 
		|F(z)|^2\leq K_E(z,z)\cdot\frac 1\pi\int\limits_{\bb R}\Big|\frac{F(t)}{E(t)}\Big|^2\DD t
		.
	\]
	If $F\in\mc H(E)$, then $\|F\|^2=\frac 1\pi\int_{\bb R}\big|\frac{F(t)}{E(t)}\big|^2\DD t$. 
\end{theorem}

\noindent
Note here that $K_E(z,z)^{\frac 12}$ is the norm of the point evaluation functional at $z$. 

The freedom of choice of $E$ expressed by \Cref{M99}\,$(iii)$ can be used to impose certain normalisations. For example, it is
always possible to choose $E$ with $E(0)=1$. 

When working with Hermite-Biehler functions, the following notion is often used. 

\begin{remark}
\label{M81}
	For a Hermite-Biehler function $E$ there exists a continuously differentiable function $\varphi_E$ with $\varphi_E'>0$, 
	such that $E(t)e^{i\varphi_E(t)}\in\bb R$ for all $t\in\bb R$. Such a function is unique up to an additive constant in
	$\pi\bb Z$, and each such function is called \emph{phase function} of $E$. 
\end{remark}

\begin{center}
	\textnumero 2
	\\
	\emph{The chain of de~Branges subspaces}
\end{center}
A central role in the theory of de~Branges spaces is played by those subspaces of a given space which are themselves de~Branges 
spaces. 

\begin{definition}
\label{M49}
	Let $\mc H$ be a de~Branges space. A linear subspace $\mc L$ of $\mc H$ is called a \emph{de~Branges subspace}, 
	if it is closed, invariant under the involution $F\mapsto F^\#$, and invariant under division of zeroes. 

	We denote the set of all de~Branges subspaces of $\mc H$ as $\Sub\mc H$. 
\end{definition}

\noindent
Note that $\mc L$ is a de~Branges subspace of $\mc H$ if and only if it is with the inner product inherited from 
$\mc H$ itself a de~Branges space. It is a significant result that $\Sub\mc H$ has a very particular order structure. 

\begin{theorem}
\label{M27}
	Let $\mc H$ be a de~Branges space. Then $\Sub\mc H$ is totally ordered with respect to inclusion. 
	For every nonempty subset $S\subseteq\Sub\mc H$ we have 
	\[
		\Cls\bigcup S\in\Sub\mc H,\quad \bigcap S\in\Sub\mc H\cup\big\{\{0\}\big\}
		,
	\]
	and it holds that 
	\[
		\forall \mc L\in\Sub\mc H\setminus\{\mc H\}\DP
		\dim\bigg[
		\raisebox{2pt}{$\displaystyle\bigcap\big\{\mc L'\DS\mc L'\in\Sub\mc H,\mc L'\supsetneq\mc L\big\}$}
		\Big/
		\raisebox{-4pt}{$\mc L$}
		\bigg]\leq 1
		,
	\]
	\[
		\forall \mc L\in\Sub\mc H,\dim\mc L>1\DP
		\dim\bigg[
		\raisebox{2pt}{$\mc L$}
		\Big/
		\raisebox{-4pt}{$\displaystyle\Cls\bigcup\big\{\mc L'\DS\mc L'\in\Sub\mc H,\mc L'\subsetneq\mc L\big\}$}
		\bigg]\leq 1
		.
	\]
\end{theorem}

\noindent
The fact that $\Sub\mc H$ is a chain is known as \emph{de~Branges' ordering theorem}, and its proof heavily relies on function
theoretic tools. 

In the context of de~Branges subspaces we can explain one important aspect of the fourth condition in \Cref{M101}. Without
imposing this condition, the set $\Sub\mc H$ is not a chain but a disjoint union of chains (provided that $\dim\mc H=\infty$,
for each possible distribution of real zeroes exactly one). Thus, requiring the fourth condition means in some sense to impose a
canonical way to pick one chain of subspaces. 

For later reference we include another property of $\Sub\mc H$. Here we denote, for any de~Branges space $\mc H$, its
reproducing kernel by $K_{\mc H}$. 

\begin{lemma}
\label{M86}
	Let $\mc H$ be a de~Branges space. Then the set 
	\[
		S\DE\big\{K_{\mc L}\DS\mc L\in\Sub\mc H\big\}\cup\{0\}
	\]
	is closed in the topology of locally uniform convergence on $\bb C^2$.
\end{lemma}
\noindent
The set $S$ is even closed w.r.t.\ pointwise convergence. This follows using a compactness argument, but 
is not needed at present.

\begin{proof}[Proof of \Cref{M86}]
	Since locally uniform convergence is metrisable, it is enough to consider limits of sequences of elements of $S$. 
	If the limit is $0$, there is nothing to prove. Hence, it is enough to consider a sequence of spaces 
	$\mc L_n\in\Sub\mc H$ such that the limit $K\DE\lim_{n\to\infty}K_{\mc L_n}$ exists locally uniformly and $K\neq 0$. 
	Since $\Sub\mc H$ is totally ordered, there exists a subsequence of $(\mc L_n)_{n=1}^\infty$ that is monotone. Hence, we
	may assume w.l.o.g.\ that $(\mc L_n)_{n=1}^\infty$ is monotone. 

	In the following we denote by $P_n$ the orthogonal of $\mc H$ onto $\mc L_n$. Then 
	\[
		\forall w\in\bb C\DP K_{\mc L_n}(\Dummy,w)=P_nK_{\mc H}(\Dummy,w).
	\]
	Consider the case that $(\mc L_n)_{n=1}^\infty$ is nondecreasing.
	Denote $\mc L_\infty\DE\Cls\bigcup_{n=1}^\infty\mc L_n$, and let $P_\infty$ be the orthogonal projection of $\mc H$ 
	onto $\mc L_\infty$. Then $\lim_{n\to\infty}P_n=P_\infty$ in the strong operator topology. Thus, 
	\begin{align*}
		K(z,w)= &\, \lim_{n\to\infty}K_{\mc L_n}(z,w)
		=\lim_{n\to\infty}\big(P_nK_{\mc H}(\Dummy,w),K_{\mc H}(\Dummy,w)\big)
		\\
		= &\, \lim_{n\to\infty}\big(P_\infty K_{\mc H}(\Dummy,w),K_{\mc H}(\Dummy,w)\big)
		=K_{\mc L_\infty}(z,w)\in S.
	\end{align*}
	The case that $(\mc L_n)_{n=1}^\infty$ is nonincreasing is treated in the same; just use 
	$\mc L_\infty\DE\bigcap_{n=1}^\infty\mc L_n$ instead.
\end{proof}

\begin{center}
	\textnumero 3
	\\
	\emph{The structure Hamiltonian}
\end{center}
Let $\mc H$ be a de~Branges space. When passing to the concrete description via Hermite-Biehler functions 
the chain $\Sub\mc H$ can be described by means of a differential equation. 

To explain this, we need to make a small excursion to the theory of canoncial systems. 
A two-dimensional \emph{canonical system} is a differential equation of the form (for practical reasons we write the 
equation for row vectors) 
\begin{equation}
\label{M102}
	\frac{\partial}{\partial t}(y_1(t),y_2(t))J=z(y_1(t),y_2(t))H(t),\qquad t\in(l_-,l_+)\text{ a.e.},
\end{equation}
where $-\infty\leq l_-<l_+\leq\infty$, $J$ is the symplectic matrix $J\DE\smmatrix 0{-1}10$, 
$z\in\bb C$ is the eigenvalue parameter, and where $H\in L^1_{\textrm{loc}}((l_-,l_+),\bb R^{2\times 2})$.
The function $H$ is called the \emph{Hamiltonian} of the system\footnote{%
	We deliberately do not assume that $H$ is positive semidefinite or that $H$ is integrable up to one or both endpoints 
	of the interval.}.

The equation \cref{M102} can be rewritten in integral form. In fact, a function 
$(y_1,y_2)\DF (l_-,l_+)\to\bb C^{1\times 2}$ is locally absolutely continuous and satisfies \cref{M102}, 
if and only if it is measurable, locally bounded, and satisfies 
\[
	\forall l_-<a<b<l_+\DP
	(y_1(b),y_2(b))J-(y_1(a),y_2(a))J=z\int\limits_a^b (y_1(t),y_2(t))H(t)\DD t
	.
\]
An interval $(a,b)\subseteq(l_-,l_+)$ is called \emph{indivisible of type $\phi$}, if 
\[
	\ker H(t)=\Span\Big\{\binom{-\sin\phi}{\cos\phi}\Big\},\quad t\in(a,b)\text{ a.e}
	.
\]
Given a Hamiltonian $H$ on an interval $(l_-,l_+)$ and a point $t\in(l_-,l_+)$, we denote by $W_H(t,s,z)$, $s\in(l_-,l_+)$, 
the unique $2\times 2$-matrix solution of the initial value problem 
\[
	\left\{
	\begin{array}{l}
		\frac{\partial}{\partial s}W_H(t,s,z)J=zW_H(t,s,z)H(s),\quad 
		s\in(l_-,l_+)\text{ a.e.},
		\\[2mm]
		W_H(t,t,z)=I.
	\end{array}
	\right.
\]
We refer to $W_H(t,s,z)$ as the family of \emph{transfer matrices} associated with $H$. 
Note that 
\begin{align*}
	& \forall t,s,r\in(l_-,l_+)\DP W_H(t,s,z)W_H(s,r,z)=W_H(t,r,z),
	\\
	& \forall t,s\in(l_-,l_+)\DP W_H(t,s,0)=I
	.
\end{align*}

\begin{theorem}
\label{M25}
	Let $E\in\HB$ with $E(0)=1$. Then there exists a unique Hamiltonian $H_E$ defined on the interval 
	$(-\infty,0)$ with
	\[
		H_E(t)\geq 0\text{ a.e.},\quad \Tr H_E(t)=1\text{ a.e.},\quad 
		\int\limits_{-\infty}^0\binom 10^*H_E(t)\binom 10\DD t<\infty
		,
	\]
	such that the solution of the initial value problem 
	\[
		\left\{
		\begin{array}{l}
			\frac{\partial}{\partial t}(A(t,z),B(t,z))J=z(A(t,z),B(t,z))H_E(t),\quad 
			t\in(-\infty,0)\text{ a.e.},
			\\[2mm]
			(A(0,z),B(0,z))=(A(z),B(z)),
		\end{array}
		\right.
	\]
	satisfies (following the generic notation \cref{M50} we write $E(t,.)\DE A(t,.)-iB(t,.)$) 
	\[
		\forall t\in(-\infty,0]\DP E(t,.)\in\HB\cup\{1\}
		,
	\]
	\phantom{}\\[-13mm]
	\begin{align*}
		\Sub\mc H(E)=\big\{\mc H(E(t,.))\DS & t\in(-\infty,0],E(t,.)\neq 1,
		\\
		& t\text{ is not inner point of an indivisible interval\,}\big\}
		.
	\end{align*}
\end{theorem}

\noindent
The Hamiltonian $H_E$, granted uniquely by this theorem, is called the \emph{structure Hamiltonian} associated with $E$, and 
we write the transfer matrices of $H_E$ as $W_E(t,s,z)$.

The freedom in the choice of $E$ when representing a de~Branges space $\mc H$ as $\mc H(E)$ expressed by \Cref{M99}\,$(iii)$ 
translates easily to structure Hamiltonians. Given the normalisation that $E(0)=1$, each two Hermite-Biehler functions 
$E,\tilde E$ with $\mc H(E)=\mc H(\tilde E)$ are related as 
\[
	(A,B)=(\tilde A,\tilde B)\smmatrix 10{\gamma}1
\]
with some real constant $\gamma$. The corresponding structure Hamiltonians are then related as 
\[
	H_E=\smmatrix 10{\gamma}1H_{\tilde E}\smmatrix 1{\gamma}01
	.
\]

\begin{center}
	\textnumero 4
	\\
	\emph{Unbounded chains and measures}
\end{center}
Given a de~Branges space $\mc H$, there always exist positive Borel measures on the real line, such that $\mc H$ is contained
isometrically in $L^2(\mu)$\footnote{%
	To make it explicit: we say that $\mc H\subseteq L^2(\mu)$ isometrically, if 
	\[
		\forall F\in\mc H\DP \|F\|^2=\int\limits_{\bb R}\big|F(t)\big|^2\DD\mu(t)
		.
	\]
	}. 
For example, choose $E\in\HB$ with $\mc H=\mc H(E)$, then 
$\mc H\subseteq L^2(\frac 1\pi\frac{\DD t}{|E(t)|^2})$ isometrically. 
Other examples are obtained using orthonormal bases of $\mc H$, and such can also be constructed explicitly from $E$. 

In the description of all measures $\mu$ such that $\mc H$ is contained isometrically in $L^2(\mu)$, chains of de~Branges spaces
occur which, unlike $\Sub\mc H$, do not have a maximal element. In the present context, the theory of such chains is not needed
in its full generality; the reason being that in the context of homogeneous spaces the dimensions occuring in \Cref{M27} 
are always equal to $0$. In the following discussion we restrict to what is needed at present.

\begin{definition}
\label{M22}
	We call a set $\mc C$ of de~Branges spaces an \emph{unbounded chain}, if 
	\begin{Itemize}
	\item $\mc C$ is totally ordered with respect to isometric inclusion;
	\item $\forall \mc H\in\mc C\DP \Sub\mc H\subseteq\mc C$;
	\item $\lim\limits_{\mc H\in\mc C}K_{\mc H}(0,0)=\infty$, where $K_{\mc H}$ denotes the reproducing kernel of $\mc H$, 
		and $\mc C$ is understood as a directed set w.r.t.\ inclusion.
	\end{Itemize}
\end{definition}

\begin{theorem}
\label{M23}
	Let $\mc C$ be an unbounded chain of de~Branges spaces. Then the following statements hold.
	\begin{Enumerate}
	\item There exists a unique positive Borel measure $\mu_{\mc C}$ on $\bb R$, such that 
		\begin{align*}
			& \forall \mc H\in\mc C\DP \mc H\subseteq L^2(\mu)\text{ isometrically}
			,
			\\
			& \bigcup\big\{\mc H\DS \mc H\in\mc C\big\}\text{ is dense in }L^2(\mu)
			.
		\end{align*}
	\item If $(\mu_{\mc H})_{\mc H\in\mc C}$ is a net of positive Borel measures such that $\mc H\subseteq L^2(\mu_{\mc H})$ 
		isometrically for all $\mc H\in\mc C$, then $\lim_{\mc H\in\mc C}\mu_{\mc H}=\mu_{\mc C}$ in the sense of vague 
		convergence of measures\footnote{%
			For this notion of convergence see, e.g., \cite[\S 13.2]{klenke:2020}.}.
	\end{Enumerate}
\end{theorem}

\noindent
Passing to the concrete level of descriptions via Hermite-Biehler functions, also unbounded chains $\mc C$ and the measure 
$\mu_{\mc C}$ can be understood with the help of canoncial systems. To formulate this fact, we need to recall the notions of
Nevanlinna functions and the Weyl coefficient of a limit point system. 

A function $q$ is called a \emph{Nevanlinna function}, if it is analytic on $\bb C\setminus\bb R$, satisfies 
$\Im q(z)\geq 0$ for all $z\in\bb C_+$, and $q(\ov z)=\ov{q(z)}$ for all $z\in\bb C\setminus\bb R$. 
Assume now we have a Hamiltonian on some interval $(l_-,l_+)$ with $H(t)\geq 0$ a.e., that is integrable at $l_-$ 
but not integrable on the whole interval. 
Then for every family $(q_t)_{t\in[l_-,l_+)}$ of Nevanlinna functions the limit\footnote{%
	Due to the assumption that $H$ is integrable at $l_-$, the transfer matrix exists also starting with inital node
	$l_-$.}
	\footnote{The symbol ``$\star$'' denotes the usual action of $\GL(2,\bb C)$ on the Riemann sphere $\bb C_\infty$:
	\[
		\smmatrix{m_{11}}{m_{12}}{m_{21}}{m_{22}}\star\tau\DE\frac{m_{11}\tau+m_{12}}{m_{21}\tau+m_{22}}
		.
	\]
	}
\[
	q_H(z)\DE\lim_{t\uparrow l_+}W_H(l_-,t,z)\star q_t
\]
exists locally uniformly on $\bb C\setminus\bb R$, is independent of the choice of $q_t$, and is a Nevanlinna function. 
This limit is called the \emph{Weyl coefficient} of $H$. 

\begin{theorem}
\label{M34}
	Let $\mc C$ be an unbounded chain of de~Branges spaces, and let $E\in\HB$ with $E(0)=1$ and $\mc H(E)\in\mc C$. 
	Then the following statements hold.
	\begin{Enumerate}
	\item There exists a unique Hamiltonian $H$ on $(0,\infty)$ with 
		\[
			H(t)\geq 0\text{ a.e.},\quad \Tr H(t)=1\text{ a.e.}
			,
		\]
		such that the solution of the initial value problem (note that by the above properties $H$ is integrable at $0$)
		\[
			\left\{
			\begin{array}{l}
				\frac{\partial}{\partial t}(A(t,z),B(t,z))J=z(A(t,z),B(t,z))H(t),\quad 
				t\in(0,\infty)\text{ a.e.},
				\\[2mm]
				(A(0,z),B(0,z))=(A(z),B(z)),
			\end{array}
			\right.
		\]
		satisfies 
		\[
			\forall t\in[0,\infty)\DP E(t,.)\in\HB
			,
		\]
		\phantom{}\\[-13mm]
		\begin{multline*}
			\big\{\mc L\in\mc C\DS\mc H(E)\subseteq\mc L\big\}=\big\{\mc H(E(t,.))\DS t\in[0,\infty),
			\\
			t\text{ is not inner point of an indivisible interval\,}\big\}
			.
		\end{multline*}
	\item Let $q_H$ be the Weyl coefficient of $H$, and set 
		\begin{equation}
		\label{M67}
			q_{E,\mc C}\DE\begin{pmatrix} A & B \\ -B & A\end{pmatrix}\star\Big(\frac{-1}{q_H}\Big)
			.
		\end{equation}
		Then there exists $\beta\geq 0$ such that 
		\[
			\Im q_{E,\mc C}(x+iy)=
			\beta y+\int\limits_{\bb R}\frac y{(t-x)^2+y^2}\cdot|E(t)|^2\DD\mu_{\mc C}(t)
			,\quad x\in\bb R,y>0
			.
		\]
		If $\bigcup\big\{\mc L\in\Sub\mc H(E)\DS\mc L\subsetneq\mc H(E)\big\}$ is dense in $\mc H(E)$, then $\beta=0$. 
	\end{Enumerate}
\end{theorem}

\noindent
The representation of $\mu_{\mc C}$ given in \Cref{M34}\,$(ii)$ yields a continuity property\footnote{%
	We denote $\bb N\DE\{1,2,3,\ldots\}$ and $\bb N_0\DE\bb N\cup\{0\}$
	}.

\begin{lemma}
\label{M62}
	Let $\mc C_n$, $n\in\bb N_0\cup\{\infty\}$, be unbounded chains.  
	Let $E_n\in\HB$, $n\in\bb N_0\cup\{\infty\}$ be such that $E_n(0)=1$ and $\mc H(E_n)\in\mc C_n$ for all 
	$n\in\bb N_0\cup\{\infty\}$, and denote by $H_n$ the corresponding Hamiltonians granted by \Cref{M34}\,$(i)$.

	Assume that $\lim_{n\to\infty}E_n=E_\infty$ locally uniformly on $\bb C$, and 
	$\lim_{n\to\infty}H_n=H_\infty$ locally weak-$L^1$ on $[0,\infty)$. Then 
	$\lim_{n\to\infty}\mu_{\mc C_n}=\mu_{\mc C_\infty}$ vaguely. 
\end{lemma}

\noindent
Since results of this kind are not discussed in \cite{debranges:1968}, we provide the argument.

\begin{proof}
	Consider the respective functions $q_{E_n,\mc C_n}$, $n\in\bb N_0\cup\{\infty\}$ introduced in \cref{M67}. 
	By our assumptions on convergence of $E_n$ and $H_n$ we have 
	$\lim_{n\to\infty}q_{E_n,\mc C_n}=q_{E_\infty,\mc C_\infty}$ locally uniformly on $\bb C\setminus\bb R$ 
	(recall here that, by a theorem about canonical systems, convergence of Hamiltonians implies convergence of 
	Weyl coefficients, e.g.\ \cite{remling:2018}). 
	The Grommer-Hamburger theorem implies that $\lim_{n\to\infty}|E_n(t)|^2\DD\mu_{\mc C_n}=
	|E_\infty(t)|^2\DD\mu_{\mc C_\infty}$ vaguely. Since the functions $E_n$ are continuous and have no real zeroes, it
	follows that $\lim_{n\to\infty}\mu_{\mc C_n}=\mu_{\mc C_\infty}$ vaguely. 
\end{proof}

%\newpage
%
%
%
\section[\textcolor{ForestGreen}{Introduction to the main players}]{Introduction to the main players}
%ForestGreen / Dandelion / BrickRed / Fuchsia
%
%
\label{M103}

\subsection[\textcolor{ForestGreen}{Weighted rescaling}]{Weighted rescaling}
%ForestGreen / Dandelion / BrickRed / Fuchsia
%

For each $p\in\bb R$ we define a continuous group action $\odot_p$ of the positive real numbers on a function space 
$X^{\bb C}$, where at first $X$ is just any normed space; later on it will mainly be $\bb C$ or $\bb C^2$. 
Here (and always) $X^{\bb C}$ is endowed with the topology of locally uniform convergence.

\begin{definition}
\label{M1}
	Let $p\in\bb R$ and $X$ a normed space. Then $\odot_p:\bb R^+\times X^{\bb C}\to X^{\bb C}$ is defined as
	\begin{equation}
	\label{M2}
		[a\odot_p F](z)\DE a^p F(az)\quad\text{for }z\in\bb C
		.
	\end{equation}
\end{definition}

\noindent
It is obvious that $\odot_p$ indeed is a continuous group action. 

\begin{Remark}
\label{M109}
	When speaking of weighted rescalings one could also think of using other weights $k(a)$ than powers 
	in \cref{M100} and \cref{M2}. For the following two reasons this does in the present context not lead to greater
	generality. 
	\begin{Enumerate}
	\item If we want to get a group action, the weight function 
		$k$ must be a solution of the multiplicative Cauchy functional equation $k(ab)=k(a)k(b)$. 
	\item If for all $a\in(0,1]$ the map $F(z)\mapsto k(a)F(az)$ maps some de~Branges space isometrically into itself, then 
		$k$ is a solution of the multiplicative Cauchy functional equation. This is proven in \cite{debranges:1962b}.
	\end{Enumerate}
	Making some weak assumption on $k$, for instance that $k$ is measurable, it is thus no loss of generality to restrict
	attention to weights $k(a)=a^p$ where $p\in\bb R$. 
\end{Remark}

The group action $\odot_p$ fits well with the construction of de~Branges spaces from Hermite-Biehler functions. 

\begin{lemma}
\label{M3}
	Let $p\in\bb R$ and $E\in\HB$, and let further $a\in\bb R^+$. 
	Then $a\odot_p E\in\HB$, and the reproducing kernels of $\mc H(E)$ and $\mc H(a\odot_p E)$ are related as 
	\begin{equation}
	\label{M4}
		K_{a\odot_p E}(z,w)=a^{2p+1}K_E(az,aw)\quad\text{for }z,w\in\bb C
		.
	\end{equation}
	The map $F\mapsto a\odot_{p+\frac 12} F$ is an isometric isomorphism of $\mc H(E)$ onto $\mc H(a\odot_p E)$. 
\end{lemma}
\begin{proof}
	The fact that $a\odot_p E$ is a Hermite-Biehler function is obvious. The decomposition of $a\odot_p E$ in real- and
	imaginary components is 
	\[
		a\odot_p E=(a\odot_p A)-i(a\odot_p B)
		,
	\]
	and from this we immediately obtain the kernel relation \cref{M4}.

	Considering $w\in\bb C$ as a fixed parameter, \cref{M4} says that 
	\[
		\big[a\odot_{p+\frac 12}K_E(.,w)\big](z)=a^{p+\frac 12}K_E(az,w)=
		a^{-p-\frac 12}K_{a\odot_p E}\big(z,\frac wa\big)
		\quad\text{for }z\in\bb C
		.
	\]
	Thus $a\odot_{p+\frac 12}K_E(.,w)\in\mc H(a\odot_p E)$, and for each two points $w,w'\in\bb C$, 
	\begin{align*}
		\Big(a\odot_{p+\frac 12} &\, K_E(.,w),a\odot_{p+\frac 12}K_E(.,w')\Big)_{\mc H(a\odot_p E)}
		\\
		= &\, a^{-2p-1}
		\Big(K_{a\odot_p E}\big(.,\frac wa\big),K_{a\odot_p E}\big(.,\frac{w'}a\big)\Big)_{\mc H(a\odot_p E)}
		\\
		= &\,a^{-2p-1}K_{a\odot_p E}\big(\frac{w'}a,\frac wa\big)
		=K_E(w',w)=\big(K_E(.,w),K_E(.,w')\big)_{\mc H(E)}
		.
	\end{align*}
	We see that $F\mapsto a\odot_{p+\frac 12}F$ maps the linear span of reproducing kernels of $\mc H(E)$ isometrically 
	onto the linear span of reproducing kernels of $\mc H(a\odot_p E)$. Hence, it extends to an isometric isomorphism 
	of $\mc H(E)$ onto $\mc H(a\odot_p E)$. Since point evaluations are continuous in both spaces, 
	this extension acts again as $F\mapsto a\odot_{p+\frac 12}F$. 
\end{proof}

\subsection[\textcolor{ForestGreen}{The canonical system}]{The canonical system}
%ForestGreen / Dandelion / BrickRed / Fuchsia
%

We define a class of Hamiltonians having a very particular form. 

\begin{Definition}
\label{M52}
	Let $p\in\bb R$ and $(P,\psi)\in\bb R^{2\times 2}\times\bb R$. Then we define functions 
	$\ms D_\psi,H_{P,\psi}\DF (0,\infty)\to\bb R^{2\times 2}$ as 
	\begin{align*}
		\ms D_\psi(a)\DE &\, 
		\begin{cases}
			\begin{pmatrix} 1 & 0 \\ \frac{\psi}{2p} & 1 \end{pmatrix}
			\begin{pmatrix} a^p & 0 \\ 0 & a^{-p} \end{pmatrix}
			\begin{pmatrix} 1 & 0 \\ -\frac{\psi}{2p} & 1 \end{pmatrix}
			\CAS p\neq 0
			,
			\\[6mm]
			\begin{pmatrix} 1 & 0 \\ \psi\log a & 1 \end{pmatrix}
			\CAS p=0
			,
		\end{cases}
		\\[3mm]
		H_{P,\psi}(a)\DE &\, \ms D_\psi(a)P\ms D_\psi(a)^T
		.
	\end{align*}
	Here we slightly overload notation by not explicitly denoting dependence on $p$.
\end{Definition}

\noindent
For each fixed $P$ and $\psi$ both of $\ms D_\psi$ and $H_{P,\psi}$ are continuous (in fact, infinitely differentiable)
in $a\in(0,\infty)$. The functions $\ms D_\psi$ (for fixed $p$) form a semigroup:
\[
	\ms D_\psi(1)=I,\qquad\forall a,a'\in(0,\infty)\DP \ms D_\psi(a)\ms D_\psi(a')=\ms D_\psi(aa').
\]
Moreover, $\ms D_\psi(a)\in\SL(2,\bb R)$ for all $a\in(0,\infty)$.
For later reference we also state explicitly that, for any $p\neq 0$, 
\begin{equation}
\label{M88}
	\ms D_\psi(a)^T\binom 10=a^p\binom 10,\quad\ms D_\psi(a)^T\binom{-\psi}{2p}=a^{-p}\binom{-\psi}{2p}.
\end{equation}
Let us have a closer look at the function $H_{P,\psi}$. First, we consider $\ker H_{P,\psi}(a)$. It is clear that 
$\Rank H_{P,\psi}(a)=\Rank P$ for all $a>0$. In particular, if $H_{P,\psi}(a)=0$ or $H_{P,\psi}(a)$ is invertible for 
one $a>0$, then the respective property holds for all $a>0$. The behaviour of $\ker H_{P,\psi}(a)$ when $\Rank P=1$ 
is slightly more complex.

\begin{lemma}
\label{M41}
	Let $p\in\bb R$ and $(P,\psi)\in\bb R^{2\times 2}\times\bb R$, and assume that $\ker P=\Span\{\xi\}$ with some 
	nonzero vector $\xi$. 
	\begin{Enumerate}
	\item Assume that $p\neq 0$. If $\xi$ is a scalar multiple of either $\binom 10$ or $\binom{-\psi}{2p}$, then 
		$\ker H_{P,\psi}(a)=\Span\{\xi\}$ for all $a>0$. Otherwise, 
		$\ker H_{P,\psi}(a)\neq\ker H_{P,\psi}(b)$ for any two $a,b$ with $0<a<b<\infty$.
	\item Assume that $p=0$ and $\psi\neq 0$. If $\xi$ is a scalar multiple of $\binom 10$, then 
		$\ker H_{P,\psi}(a)=\Span\{\xi\}$ for all $a>0$. Otherwise, 
		$\ker H_{P,\psi}(a)\neq\ker H_{P,\psi}(b)$ for any two $a,b$ with $0<a<b<\infty$.
	\item Assume that $p=\psi=0$. Then $\ker H_{P,\psi}(a)=\Span\{\xi\}$ for all $a>0$. 
	\end{Enumerate}
\end{lemma}
\begin{proof}
	Since $\ms D_\psi(a)$ is invertible, we have 
	\[
		\ker H_{P,\psi}(a)=\Span\big\{\ms D_\psi(a)^{-T}\xi\big\}
		.
	\]
	Consider the case that $p\neq 0$. Then 
	\[
		\ms D_\psi(a)^{-T}=
		\begin{pmatrix}
			1 & -\frac{\psi}{2p}
			\\
			0 & 1
		\end{pmatrix}
		\begin{pmatrix}
			a^{-p} & 0
			\\
			0 & a^p
		\end{pmatrix}
		\begin{pmatrix}
			1 & \frac{\psi}{2p}
			\\
			0 & 1
		\end{pmatrix}
		.
	\]
	The vectors $\binom 10$ and $\binom{-\psi}{2p}$ are eigenvectors of this matrix with respective eigenvalues 
	$a^{-p}$ and $a^p$. The first statement in (i) readily follows. Assume that $a<b$ and 
	$\ker H_{P,\psi}(a)=\ker H_{P,\psi}(b)$. Then $\ms D_\psi(a)^{-T}\xi=\alpha\ms D_\psi^{-T}(b)\xi$ for some 
	$\alpha\in\bb R\setminus\{0\}$. By multiplicativity $\xi$ is an eigenvector of $\ms D_\psi(\frac ba)^T$, and since 
	$\frac ba\neq 1$ therefore a scalar multiple of either $\binom 10$ or $\binom{-\psi}{2p}$. This finishes the proof of
	(i). 

	Consider next the case that $p=0$ and $\psi\neq 0$. Then 
	\[
		\ms D_\psi(a)^{-T}=
		\begin{pmatrix}
			1 & -\psi\log a
			\\
			0 & 1
		\end{pmatrix}
		,
	\]
	and the assertion of $(ii)$ follows immediately. 
	Also $(iii)$ is clear, since for $p=\psi=0$ we have $\ms D_\psi(a)=I$ for all $a>0$. 
\end{proof}

\noindent
Second, we investigate integrability of $H_{P,\psi}$ at the endpoints of the interval $(0,\infty)$. 
The proof is simply by explicit consideration.

\begin{lemma}
\label{M60}
	Let $(P,\psi)\in\bb R^{2\times 2}\times\bb R$. 
	Concerning integrability at $0$ we have 
	\begin{align*}
		& \forall p\in(-\infty,-\tfrac 12]\DP
		\bigg[H_{P,\psi}\in L^1\big((0,1),\bb R^{2\times 2}\big)\ \Leftrightarrow\ 
		\binom 10^*P\binom 10=0\bigg]
		,
		\\
		& \forall p\in(-\tfrac 12,\tfrac 12)\DP
		H_{P,\psi}\in L^1\big((0,1),\bb R^{2\times 2}\big)
		,
		\\
		& \forall p\in[\tfrac 12,\infty)\DP
		\bigg[H_{P,\psi}\in L^1\big((0,1),\bb R^{2\times 2}\big)\ \Leftrightarrow\ 
		\binom{-\psi}{2p}^*P\binom{-\psi}{2p}=0\bigg]
		,
	\end{align*}
	and concerning integrability at $\infty$ 
	\begin{align*}
		& \forall p\in(-\infty,-\tfrac 12)\DP
		\bigg[H_{P,\psi}\in L^1\big((1,\infty),\bb R^{2\times 2}\big)\ \Leftrightarrow\ 
		\binom{-\psi}{2p}\in\ker P\cap\ker P^T\bigg]
		,
		\\
		& \forall p\in[-\tfrac 12,\tfrac 12]\DP
		\bigg[H_{P,\psi}\in L^1\big((1,\infty),\bb R^{2\times 2}\big)\ \Leftrightarrow\ P=0\bigg]
		,
		\\
		& \forall p\in(\tfrac 12,\infty)\DP
		\bigg[H_{P,\psi}\in L^1\big((1,\infty),\bb R^{2\times 2}\big)\ \Leftrightarrow\ 
		\binom 10\in\ker P\cap\ker P^T\bigg]
		.
	\end{align*}
\end{lemma}
\begin{proof}
	Let us first settle the case that $p=0$. We write $P=\smmatrix{p_{11}}{p_{12}}{p_{21}}{p_{22}}$ and write out the
	definition of $H_{P,\psi}(a)$. This yields 
	\[
		H_{P,\psi}(a)=
		\begin{pmatrix}
			p_{11} & p_{11}\psi\log a+p_{12}
			\\[2mm]
			p_{11}\psi\log a+p_{21} & p_{11}(\psi\log a)^2+(p_{12}+p_{21})\psi\log a+p_{22}
		\end{pmatrix}
		.
	\]
	From this it is clear that $H_{P,\psi}\in L^1(0,1)$ but $H_{P,\psi}\notin L^1(1,\infty)$ unless $P=0$. 

	From now on assume that $p\neq 0$. Set $f\DE\binom 10$ and $g\DE\binom{-\psi}{2p}$. Then $H_{P,\psi}$ is integrable at
	$0$ (or at $\infty$) if and only if all four of 
	\[
		f^TH_{P,\psi}f,\ g^TH_{P,\psi}g,\ f^TH_{P,\psi}g,\ g^TH_{P,\psi}f
	\]
	have the respective property. Using \cref{M88} we obtain 
	\begin{align*}
		& f^TH_{P,\psi}f=a^{2p}\cdot f^TPf,\quad g^TH_{P,\psi}g=a^{-2p}\cdot g^TPg,
		\\
		& f^TH_{P,\psi}g=f^TPg,\quad g^TH_{P,\psi}f=g^TPf.
	\end{align*}
	This shows that 
	\[
	\begin{array}{c|c|c}
		& \parbox{30mm}{integrable\\ at $0$ $\ \Leftrightarrow$} 
		& \parbox{30mm}{integrable\\ at $\infty$ $\ \Leftrightarrow$}
		\\[3mm]
		\hline
		f^TH_{P,\psi}f\raisebox{12pt}{\phantom{.}} 
		& 2p>-1\vee f^TPf=0
		& 2p<-1\vee f^TPf=0
		\\[2mm]
		g^TH_{P,\psi}g
		& -2p>-1\vee g^TPg=0
		& -2p<-1\vee g^TPg=0
		\\[2mm]
		f^TH_{P,\psi}g
		& \text{always}
		& f^TPg=0
		\\[2mm]
		g^TH_{P,\psi}f
		& \text{always}
		& g^TPf=0
		\\[1mm]
		\hline
	\end{array}
	\]
	To finish the proof we simply go through all cases. For integrability at $0$:
	\begin{Itemize}
	\item $p\in(-\infty,-\tfrac 12]$: 
		Then ``$2p>-1$'' fails and ``$-2p>-1$'' holds. 
		Hence, $H_{P,\psi}$ is integrable at $0$ if and only if $f^TPf=0$. 
	\item $p\in(-\tfrac 12,\tfrac 12)$: 
		Then ``$2p>-1$'' and ``$-2p>-1$'' both hold. 
		Hence, $H_{P,\psi}$ is integrable at $0$. 
	\item $p\in[\tfrac 12,\infty)$: 
		Then ``$2p>-1$'' holds and ``$-2p>-1$'' fails. 
		Hence, $H_{P,\psi}$ is integrable at $0$ if and only if $g^TPg=0$. 
	\end{Itemize}
	For integrability at $\infty$:
	\begin{Itemize}
	\item $p\in(-\infty,-\tfrac 12)$: 
		Then ``$2p<-1$'' holds and ``$-2p<-1$'' fails. 
		Hence, 
		\[
			H_{P,\psi}\text{ integrable at $\infty$}\ \Leftrightarrow\ g^TPg=f^TPg=g^TPf=0
			\Leftrightarrow\ Pg=g^TP=0
			.
		\]
	\item $p\in[-\tfrac 12,\tfrac 12]$: 
		Then ``$2p<-1$'' and ``$-2p<-1$'' both fail. 
		Hence, $H_{P,\psi}$ is integrable at $\infty$ if and only if $P=0$. 
	\item $p\in(\tfrac 12,\infty)$: 
		Then ``$2p<-1$'' fails and ``$-2p<-1$'' holds. 
		Hence, 
		\[
			H_{P,\psi}\text{ integrable at $\infty$}\ \Leftrightarrow\ f^TPf=f^TPg=g^TPf=0
			\Leftrightarrow\ Pf=f^TP=0
			.
		\]
	\end{Itemize}
\end{proof}

\subsection[\textcolor{ForestGreen}{The recurrence relation}]{The recurrence relation}
%ForestGreen / Dandelion / BrickRed / Fuchsia
%

Given parameters $p\in\bb R\setminus(-\frac 12\bb N)$ and $(P,\psi)\in\bb R^{2\times 2}\times\bb R$, we consider the
linear recurrence for a sequence $((\alpha_n,\beta_n))_{n\in\bb N_0}$ of pairs of real numbers given by 
\begin{equation}
\label{M30}
	\left\{
	\begin{array}{l}
		(\alpha_{n+1},\beta_{n+1})=(\alpha_n,\beta_n)\cdot\frac{-1}{(n+1)(2p+n+1)}PJ\smmatrix{2p+n+1}{0}{\psi}{n+1}
		\quad\text{for }n\in\bb N_0
		,
		\\[4mm]
		(\alpha_0,\beta_0)=(1,0)
		,
	\end{array}
	\right.
\end{equation}
where again $J=\smmatrix 0{-1}10$.

The solution $((\alpha_n,\beta_n))_{n\in\bb N_0}$ can be estimated. 
We have\footnote{%
	We always use the euclidean norm on $\bb C^d$, and the
	corresponding operator norm for matrices.}
\[
	C\DE \sup_{n\in\bb N_0}\bigg\|
	\begin{pmatrix}
		1 & 0
		\\
		\frac{\psi}{2p+n+1} & \frac{n+1}{2p+n+1}
	\end{pmatrix}
	\bigg\|<\infty
	,
\]
and hence obtain inductively that 
\begin{equation}
\label{M35}
	\forall n\in\bb N_0\DP \big\|(\alpha_n,\beta_n)\big\|\leq\frac{C^n\|P\|^n}{n!}
	.
\end{equation}
This shows that the generating functions 
\[
	A(z)\DE\sum_{n=0}^\infty \alpha_nz^n,\quad B(z)\DE\sum_{n=0}^\infty \beta_nz^n
	,
\]
of the sequences $(\alpha_n)_{n\in\bb N_0}$ and $(\beta_n)_{n\in\bb N_0}$ are entire and that 
\[
	|A(z)|,|B(z)|\leq e^{C\|P\|\cdot|z|}\quad\text{for }z\in\bb C
	,
\]
i.e., $A$ and $B$ are of finite exponential type not exceeding $C\|P\|$. 

In the following we denote by $\Hol(\bb C)$ the set of all entire functions. 

\begin{Definition}
\label{M83}
	Let $p\in\bb R\setminus(-\frac 12\bb N)$. We define maps
	\begin{align*}
		\Xi_p\DF &\, \bb R^{2\times 2}\times\bb R\to \Hol(\bb C)\times\Hol(\bb C)
		\\
		\widehat\Xi_p\DF &\, \bb R^{2\times 2}\times\bb R\to\Hol(\bb C)
	\end{align*}
	by assigning to a parameter $(P,\psi)$ the pair $(A,B)$ of generating functions of the solution of \cref{M30}, and 
	setting $\widehat\Xi_p(P,\psi)\DE A-iB$.
\end{Definition}

\noindent
We start with a simple but practical observation. Here (and always) we topologise domain and codomain of $\Xi_p$ and 
$\widehat\Xi_p$ naturally with the euclidean topology and the topology of locally uniform convergence, respectively. 

\begin{lemma}
\label{M63}
	The map 
	\[
		\left\{
		\begin{array}{rcl}
			\bb R\setminus(-\smfrac 12\bb N)\times\bb R^{2\times 2}\times\bb R & \to & 
			\Hol(\bb C)\times\Hol(\bb C)
			\\
			(p,P,\psi) & \mapsto & \Xi_p(P,\psi)
		\end{array}
		\right.
	\]
	is continuous. 
\end{lemma}
\begin{proof}
	Assume we have $((p_i,P_i,\psi_i))_{i\in\bb N_0}$ and $(p,P,\psi)$ with $\lim_{i\to\infty}(p_i,P_i,\psi_i)=(p,P,\psi)$, 
	and let $(\alpha_{i,n},\beta_{i,n})$ and $(\alpha_n,\beta_n)$ be the corresponding solutions of
	\cref{M30}. We have 
	\[
		C_1\DE \sup_{i\in\bb N_0}\|P_i\|<\infty,\quad 
		C_2\DE \sup_{i\in\bb N_0}\sup_{n\in\bb N_0}\bigg\|
		\begin{pmatrix}
			1 & 0
			\\
			\frac{\psi_i}{2p_i+n+1} & \frac{n+1}{2p_i+n+1}
		\end{pmatrix}
		\bigg\|<\infty
		,
	\]
	and get the uniform bound 
	\[
		\forall i,n\in\bb N_0\DP
		\|(\alpha_{i,n},\beta_{i,n})\|,\|(\alpha_n,\beta_n)\|\leq \frac{(C_1C_2)^n}{n!}
		.
	\]
	Clearly, $\lim_{i\to\infty}(\alpha_{i,n},\beta_{i,n})=(\alpha_n,\beta_n)$ for all $n\in\bb N_0$, and it 
	follows that $\lim_{i\to\infty}\Xi_{p_i}(P_i,\psi_i)=\Xi_p(P,\psi)$ locally uniformly. 
\end{proof}

\noindent
Let us collect some facts about the recurrence \cref{M30} and its generating functions. 
The case that the left upper corner of the parameter $P$ vanishes is exceptional.

\begin{lemma}
\label{M59}
	Let $p\in\bb R\setminus(-\frac 12\bb N)$ and $(P,\psi)\in\bb R^{2\times 2}\times\bb R$, 
	and let $((\alpha_n,\beta_n))_{n\in\bb N_0}$ be the unique solution of \cref{M30}. 
	Then the following statements are equivalent.
	\begin{Enumerate}
	\item $\binom 10^*P\binom 10=0$.
	\item $\beta_1=0$.
	\item $B=0$.
	\end{Enumerate}
\end{lemma}
\begin{proof}
	Set $\kappa_{11}\DE\binom 10^*P\binom 10$ and $\kappa_{21}\DE\binom 01^*P\binom 10$. Equating the second component of
	\cref{M30} for $n=0$ gives $\beta_1(1+2p)=\kappa_{11}$. This shows ``$(i)\Leftrightarrow(ii)$''. 
	Assume that $\kappa_{11}=0$. Multiplying \cref{M30} from the right with $\binom 01$ yields 
	\[
		\forall n\in\bb N_0\DP \beta_{n+1}=\frac{\kappa_{21}}{2p+n+1}\beta_n
		,
	\]
	and we conclude that $\beta_n=0$ for all $n\in\bb N_0$. 
\end{proof}

\noindent
It is an important property of a parameter $P$ in \cref{M30} whether or not $P$ is symmetric. 
One reason is that under the assumption of symmetry, the parameter $(P,\psi)$ can be recovered from the solution 
$((\alpha_n,\beta_n))_{n\in\bb N_0}$ of \cref{M30} by simple formulae. 

\begin{Definition}
\label{M134}
	We denote
	\[
		\bb P\DE \Big\{(P,\psi)\in\bb R^{2\times 2}\times\bb R\DS 
		P=P^T\text{ and }\,\,{\textstyle\binom 10^*P\binom 10\neq 0}\Big\}
		.
	\]
\end{Definition}

\begin{lemma}
\label{M31}
	Let $p\in\bb R\setminus(-\frac 12\bb N)$ and $(P,\psi)\in\bb P$.
	and write $P=\smmatrix{\kappa_1}{\kappa_3}{\kappa_3}{\kappa_2}$. 
	Then $(P,\psi)$ can be recovered from the first two terms $(\alpha_1,\beta_1),(\alpha_2,\beta_2)$ 
	of the solution $((\alpha_n,\beta_n))_{n\in\bb N_0}$ of \cref{M30} by the formulae
	\begin{align*}
		\kappa_1= &\, \beta_1(1+2p)
		,
		\\
		\kappa_2= &\, -\frac{\alpha_1\beta_2}{\beta_1^2}2(1+2p)+\frac{\alpha_1^2}{\beta_1}(1+2p)+
		\frac{\beta_2^2}{\beta_1^3}(2+2p)-2\frac{\alpha_2}{\beta_1}
		,
		\\
		\kappa_3= &\, \frac{\beta_2}{\beta_1}(2+2p)-\alpha_1(1+2p)
		,
		\\
		\psi= &\, \frac{\beta_2}{\beta_1^2}(2+2p)-\frac{\alpha_1}{\beta_1}2p
		.
	\end{align*}
	Note here that $\beta_1\neq 0$. 
\end{lemma}
\begin{proof}
	The relation \cref{M30} written for $n=0$ and $n=1$ gives the four equations 
	\[
		\begin{array}{l@{\qquad}l}
			\alpha_1-\beta_1\psi=-\kappa_3 & \beta_1(1+2p)=\kappa_1
			\\[2mm]
			2\alpha_2-\beta_2\psi=-\alpha_1\kappa_3-\beta_1\kappa_2 & \beta_2(2+2p)=\alpha_1\kappa_1+\beta_1\kappa_3
		\end{array}
	\]
	The second equation is the asserted formula for $\kappa_1$, and the fourth equation gives $\kappa_3$.
	Then we use the first equation to compute $\psi$, and the third for $\kappa_2$. 
\end{proof}

\noindent
This lemma has the following obvious consequence. 

\begin{Corollary}
\label{M72}
	Let $p\in\bb R\setminus(-\frac 12\bb N)$.
	Then $\Xi_p|_{\bb P}$ is a homeomorphism of $\bb P$ onto its image $\Xi_p(\bb P)$. 
	The same holds for $\widehat\Xi_p$. 
\end{Corollary}

\noindent
Next, we present a simple transformation which is often practical. 

\begin{lemma}
\label{M73}
	Let $p\in\bb R\setminus(-\frac 12\bb N)$ and $(P,\psi)\in\bb R^{2\times 2}\times\bb R$. Then, for all 
	$\gamma\in\bb R$, 
	\[
		\Xi_p(P,\psi)=
		\Xi_p\Big(\smmatrix 10{\gamma}1P\smmatrix 1{\gamma}01,\psi+2p\gamma\Big)\smmatrix 10\gamma 1
		.
	\]
\end{lemma}
\begin{proof}
	Let $((\alpha_n,\beta_n))_{n\in\bb N_0}$ be the solution of \cref{M30} for $(P,\psi)$. Then we compute 
	\begin{align*}
		{\scriptstyle -(n+1)} &\, {\scriptstyle (2p+n+1)}\cdot(\alpha_{n+1},\beta_{n+1})\smmatrix 10{-\gamma}1
		=(\alpha_n,\beta_n)PJ\smmatrix{2p+n+1}0\psi{n+1}\smmatrix 10{-\gamma}1
		\\[2mm]
		= &\, 
		(\alpha_n,\beta_n)\cdot
		\smmatrix 10{-\gamma}1\smmatrix 10\gamma 1\cdot P\cdot \smmatrix 1\gamma 01J\smmatrix 10\gamma 1\cdot
		\smmatrix{2p+n+1}0\psi{n+1}\smmatrix 10{-\gamma}1
		\\[2mm]
		= &\, 
		(\alpha_n,\beta_n)\smmatrix 10{-\gamma}1\cdot \smmatrix 10\gamma 1P\smmatrix 1\gamma 01\cdot
		J\cdot\smmatrix{2p+n+1}0{\psi+2p\gamma}{n+1}
		.
	\end{align*}
\end{proof}

\noindent
Using this transformation we can characterise a symmetry property of the generating functions. This result is of relevance in the
context of de~Branges spaces which are symmetric about the origin in the sense of \cite[Chapter~47]{debranges:1968}. 

\begin{lemma}
\label{M133}
	Let $p\in\bb R\setminus(-\frac 12\bb N)$ and $(P,\psi)\in\bb P$,
	and write $P=\smmatrix{\kappa_1}{\kappa_3}{\kappa_3}{\kappa_2}$. 
	Moreover, set 
	\[
		\sigma\DE 2p\kappa_3-\psi\kappa_1
		.
	\]
	Then the following statements are equivalent.
	\begin{Enumerate}
	\item $\sigma=0$;
	\item $B$ is odd and there exists $\gamma\in\bb R$ such that $A(z)-A(-z)=\gamma B(z)$;
	\item $B''(0)=0$.
	\end{Enumerate}
\end{lemma}
\begin{proof}
	As a preparation for the proof we study what happens when we make two steps at once in the recurrence \cref{M30}. 
	Specifically, we are interested when the property ``$\beta_n=0$'' is preserved. 
	Since $P^T=P$, the matrix $PJP$ is antisymmetric, and hence $PJPJ$ is diagonal. Thus we have, for every $n\in\bb N_0$, 
	\begin{align*}
		PJ\smmatrix{2p+n+1}{0}{\psi}{n+1} &\, PJ\smmatrix{2p+n+2}{0}{\psi}{n+2}\text{ is lower triangular}
		\\
		\Leftrightarrow\ &\, PJ\underbrace{\smmatrix{2p+n+1}{0}{\psi}{n+1}}_{=\smmatrix{2p}0\psi 0+(n+1)I}PJ
		\text{ is lower triangular}
		\\
		\Leftrightarrow\ &\, PJ\smmatrix{2p}0\psi 0PJ \text{ is lower triangular}.
	\end{align*}
	The right upper entry of the last matrix is equal to $-\kappa_1\sigma$, and hence it is lower triangular if and only if 
	$\sigma=0$. 

	We now come to the proof of the stated equivalences. Since $(\alpha_0,\beta_0)=(1,0)$, we obtain the implications 
	\[
		\Big[\sigma=0\ \Rightarrow\ \forall k\in\bb N_0\DP\beta_{2k}=0\Big]\quad\text{and}\quad
		\Big[\beta_2=0\ \Rightarrow\ \sigma=0\Big]
	\]
	Thus we have 
	\[
		\sigma=0\ \Leftrightarrow\ B\text{ is odd}\ \Leftrightarrow\ B''(0)=0.
	\]
	It remains to show that ``$\sigma=0$'' implies that $A(z)-A(-z)=\gamma B(z)$ for some $\gamma\in\bb R$. Assume that 
	$\sigma=0$. Applying \Cref{M73} with $\gamma\DE-\frac{\kappa_3}{\kappa_1}$ yields
	\[
		\Xi_p(P,\psi)=
		\underbrace{\Xi_p\Big(\smmatrix{\kappa_1}00{\frac{\det P}{\kappa_1}},0\Big)}_{\ED(\tilde A,\tilde B)}\cdot
		\smmatrix 10{-\frac{\kappa_3}{\kappa_1}}1
		.
	\]
	It is clear from \cref{M30} that in general for a diagonal matrix and parameter $0$ the solution will satisfy 
	$\alpha_n=0$ for all odd $n$ and $\beta_n=0$ for all even $n$. Thus $\tilde A$ is even and $\tilde B$ is odd. 
	The assertion follows since 
	\[
		A=\tilde A-\frac{\kappa_3}{\kappa_1}\tilde B,\quad B=\tilde B.
	\]
\end{proof}

%\newpage
%
%
%
\section[\textcolor{ForestGreen}{Solution of the canonical system via power series coefficients}]{Solution of the canonical 
	system via power series coefficients}
%ForestGreen / Dandelion / BrickRed / Fuchsia
%
%

In this section we give the connection between the group action $\odot_p$, canonical systems with Hamiltonians of the form 
$H_{P,\psi}$, and recurrences of the form \cref{M30}.
This is the approach used in \cite{debranges:1968}.

\begin{theorem}
\label{M53}
	Let $p\in\bb R\setminus(-\frac 12\bb N)$ and $(P,\psi)\in\bb R^{2\times 2}\times\bb R$, 
	and let $(\alpha_n)_{n\in\bb N_0}$ and $(\beta_n)_{n\in\bb N_0}$ be sequences of real numbers. 
	Then the following statements are equivalent.
	\begin{Enumerate}
	\item The sequence $((\alpha_n,\beta_n))_{n\in\bb N_0}$ satisfies \cref{M30}.
	\item The power series $A(z)\DE\sum_{n=0}^\infty \alpha_nz^n$ and $B(z)\DE\sum_{n=0}^\infty \beta_nz^n$ represent entire
		functions with $A(0)=1$ and $B(0)=0$. Using the notation 
		\begin{equation}
		\label{M65}
			\big(A(a,z),B(a,z)\big)\DE\big([a\odot_p A](z),[a\odot_p B](z)\big)\ms D_\psi(a)^{-1}
		\end{equation}
		for $a>0$, it holds that 
		\begin{multline}
		\label{M54}
			\forall 0<a<b<\infty\DP
			\big(A(b,z),B(b,z)\big)J-\big(A(a,z),B(a,z)\big)J
			\\
			=z\int_a^b \big(A(c,z),B(c,z)\big)H_{P,\psi}(c)\DD c
			.
		\end{multline}
	\item The power series $A(z)\DE\sum_{n=0}^\infty \alpha_nz^n$ and $B(z)\DE\sum_{n=0}^\infty \beta_nz^n$ have positive radius of
		convergence, we have $A(0)=1$ and $B(0)=0$, and there exist $a,b$ with $0<a<b<\infty$ such that the 
		equality in \cref{M54} holds. 
	\end{Enumerate}
\end{theorem}
\begin{proof}
	We are going to show that ``$(i)\Rightarrow(ii)$'' and ``$(iii)\Rightarrow(i)$''. 
	The implication ``$(ii)\Rightarrow(iii)$'' is trivial. 

	We already saw in \Cref{M103} that $(i)$ implies that $A$ and $B$ are entire functions. 
	To proceed with the proof it is convenient to rewrite \cref{M54}: 
	using the definition of $H_{P,\psi}(a)$ and the fact that 
	$\ms D_\psi(a)^TJ=J\ms D_\psi(a)^{-1}$, \cref{M54} is equivalent to 
	\begin{multline}
	\label{M55}
		\big([b\odot_p A](z),[b\odot_p B](z)\big)\ms D_\psi(b)^{-1}
		-\big([a\odot_p A](z),[a\odot_p B](z)\big)\ms D_\psi(a)^{-1}
		\\
		=-z\int_a^b \big([c\odot_p A](z),[c\odot_p B](z)\big)PJ\ms D_\psi(c)^{-1}\DD c
		.
	\end{multline}
	Plugging the power series into this relation and comparing power series coefficients yields that \cref{M55} is
	equivalent to 
	\begin{multline}
	\label{M56}
		\forall n\in\bb N_0\DP
		(\alpha_{n+1},\beta_{n+1})\big[b^{p+n+1}\ms D_\psi(b)^{-1}-a^{p+n+1}\ms D_\psi(a)^{-1}\big]
		\\
		=-(\alpha_n,\beta_n)PJ\cdot\int_a^b c^{p+n}\ms D_\psi(c)^{-1}\DD c
		.
	\end{multline}
	The square bracket on the left side of this relation computes as 
	\[
		\begin{cases}
			\begin{pmatrix}
				1 & 0 
				\\
				\frac\psi{2p} & 1
			\end{pmatrix}
			\begin{pmatrix}
				b^{n+1}-a^{n+1} & 0
				\\
				0 & b^{2p+n+1}-a^{2p+n+1}
			\end{pmatrix}
			\begin{pmatrix}
				1 & 0 
				\\
				-\frac\psi{2p} & 1
			\end{pmatrix}
			\CAS p\neq 0
			,
			\\[8mm]
			\begin{pmatrix}
				b^{n+1}-a^{n+1} & 0
				\\
				-\psi\big(b^{n+1}\log b-a^{n+1}\log a\big) & b^{n+1}-a^{n+1}
			\end{pmatrix}
			\CAS p=0
			,
		\end{cases}
	\]
	and the integral on the right side as 
	\[
		\begin{cases}
			\begin{pmatrix}
				1 & 0 
				\\
				\frac\psi{2p} & 1
			\end{pmatrix}
			\begin{pmatrix}
				\frac{b^{n+1}-a^{n+1}}{n+1} & 0
				\\
				0 & \frac{b^{2p+n+1}-a^{2p+n+1}}{2p+n+1}
			\end{pmatrix}
			\begin{pmatrix}
				1 & 0 
				\\
				-\frac\psi{2p} & 1
			\end{pmatrix}
			\CAS p\neq 0
			,
			\\[8mm]
			\begin{pmatrix}
				\frac{b^{n+1}-a^{n+1}}{n+1} & 0
				\\
				-\psi\Big(
				\frac{b^{n+1}}{n+1}\big(\log b-\frac 1{n+1}\big)-\frac{a^{n+1}}{n+1}\big(\log a-\frac 1{n+1}\big)
				\Big)
				& \frac{b^{n+1}-a^{n+1}}{n+1}
			\end{pmatrix}
			\CAS p=0
			.
		\end{cases}
	\]
	For $p\neq 0$ we can rewrite the formula for the integral as 
	\begin{multline}
	\label{M57}
		\begin{pmatrix}
			1 & 0 
			\\
			\frac\psi{2p} & 1
		\end{pmatrix}
		\begin{pmatrix}
			\frac 1{n+1} & 0
			\\
			0 & \frac 1{2p+n+1}
		\end{pmatrix}
		\begin{pmatrix}
			1 & 0 
			\\
			-\frac\psi{2p} & 1
		\end{pmatrix}
		\\
		\cdot
		\begin{pmatrix}
			1 & 0 
			\\
			\frac\psi{2p} & 1
		\end{pmatrix}
		\begin{pmatrix}
			b^{n+1}-a^{n+1} & 0
			\\
			0 & b^{2p+n+1}-a^{2p+n+1}
		\end{pmatrix}
		\begin{pmatrix}
			1 & 0 
			\\
			-\frac\psi{2p} & 1
		\end{pmatrix}
		,
	\end{multline}
	and for $p=0$ as  
	\begin{equation}
	\label{M58}
		\frac 1{n+1}
		\begin{pmatrix}
			1 & 0
			\\
			\frac\psi{n+1} & 1
		\end{pmatrix}
		\cdot
		\begin{pmatrix}
			b^{n+1}-a^{n+1} & 0
			\\
			-\psi\big(b^{n+1}\log b-a^{n+1}\log a\big) & b^{n+1}-a^{n+1}
		\end{pmatrix}
		.
	\end{equation}
	The matrix from the square bracket in \cref{M56} now appears on both sides. It is invertible since $a<b$, and 
	after cancelling out there remains a matrix on the right side which does not depend on $a$ and $b$. It equals
	\[
		\frac 1{(n+1)(2p+n+1)}
		\begin{pmatrix}
			2p+n+1 & 0
			\\
			\psi & n+1
		\end{pmatrix}
		.
	\]
\end{proof}

\noindent
From the above theorem we obtain two simple but important properties of $A$ and $B$.

\begin{corollary}
\label{M104}
	Let $p\in\bb R\setminus(-\frac 12\bb N)$ and $(P,\psi)\in\bb R^{2\times 2}\times\bb R$, and denote (as usual)
	$(A,B)\DE\Xi_p(P,\psi)$. Then $A$ and $B$ have no common zeroes. 
\end{corollary}
\begin{proof}
	The functions $A(a,z)$ and $B(a,z)$ from \cref{M65} compute explicitly as 
	\begin{align}
	\label{M106}
		A(a,z)= &\, 
		\begin{cases}
			A(az)+z\cdot\frac\psi{2p}\big(a-a^{2p+1}\big)\frac{B(az)}{az}
			\CAS p\neq 0
			,
			\\[2mm]
			A(az)-z\cdot \psi a\log a\frac{B(az)}{az}
			\CAS p=0
			,
		\end{cases}
		\\[2mm]
	\label{M107}
		B(a,z)= &\, z\cdot a^{2p+1}\frac{B(az)}{az}
		.
	\end{align}
	We see that 
	\begin{equation}
	\label{M105}
		\begin{cases}
			\lim_{a\downarrow 0}\big[A(a,z)+\frac\psi{2p}B(a,z)\big]=A(0)=1 \CAS p\neq 0
			,
			\\[1mm]
			\lim_{a\downarrow 0} A(a,z)=A(0)=1 \CAS p=0
			.
		\end{cases}
	\end{equation}
	Let $W(a,z)$ be the unique solution of the initial value problem
	\[
		\left\{
		\begin{array}{l}
			\frac{\partial}{\partial a}W(a,z)J=zW(a,z)H_{P,\psi}(a)\quad\text{for }a\in(0,\infty)
			,
			\\[2mm]
			W(1,z)=I
			.
		\end{array}
		\right.
	\]
	Note here that the initial value is prescribed at the point $a=1$, and that $A(1,z)=A(z)$ and $B(1,z)=B(z)$. 
	Then, by uniqueness of solutions, 
	\[
		\forall a\in(0,\infty)\DP
		\big(A(a,z),B(a,z)\big)=\big(A(z),B(z)\big)\cdot W(a,z)
		.
	\]
	Assume now that $z\in\bb C$ with $A(z)=B(z)=0$. Then $A(a,z)=B(a,z)=0$ for all $a>0$. 
	This contradicts \cref{M105}.
\end{proof}

\begin{corollary}
\label{M68}
	Let $p>-\frac 12$ and $(P,\psi)\in\bb R^{2\times 2}\times\bb R$, denote again $(A,B)\DE\Xi_p(P,\psi)$, 
	and let $(A(a,z),B(a,z))$ for $a>0$ be as in \cref{M65}. Then 
	\[
		\lim_{a\downarrow 0}A(a,.)=1,\quad \lim_{a\downarrow 0}B(a,.)=0
		.
	\]
	In particular, the canonical system with Hamiltonian $H_{P,\psi}$ has a solution whose limit at $0$ exists and 
	is equal to $(1,0)$. 
\end{corollary}
\begin{proof}
	This is obvious from \cref{M106} and \cref{M107}. 
\end{proof}

\begin{Remark}
\label{M108}
	In this context let us point out that any canonical system can have at most one solution 
	$(y_1(a),y_2(a))$ with $\lim_{a\downarrow l_-}(y_1(a),y_2(a))=(1,0)$. 
	This is a standard consequence of constancy of the Wronskian; for completeness we recall the argument.

	Let $-\infty\leq l_-<l_+\leq\infty$ and let $H$ be a Hamiltonian on $(l_-,l_+)$. Further, let 
	$(\eta_1,\eta_2)^T\in\bb C^2\setminus\{0\}$. Then there exists at most one solution 
	$(y_1(a),y_2(a))$ of the canonical system with Hamiltonian $H$ such that 
	$\lim_{a\downarrow l_-}(y_1(a),y_2(a))=(\eta_1,\eta_2)$. 

	To see this let $(y_1,y_2)$ and $(\tilde y_1,\tilde y_2)$ be two solutions of \cref{M102}, and assume that 
	\[
		\lim_{a\downarrow l_-}(y_1(a),y_2(a))=\lim_{a\downarrow l_-}(\tilde y_1(a),\tilde y_2(a))=(\eta_1,\eta_2)
		.
	\]
	Using the differential equation we obtain that the derivative of the Wronskian
	\[
		\det
		\begin{pmatrix}
			\tilde y_1(a) & \tilde y_2(a)
			\\
			y_1(a) & y_2(a)
		\end{pmatrix}
		=(y_1(a),y_2(a))J(\tilde y_1(a),\tilde y_2(a))^T
	\]
	is equal to zero. Hence, this determinant is independent of $a\in(l_-,l_+)$. Evaluating the limit at $l_-$ shows that 
	it is equal to $0$. Now choose $c\in(a,b)$, then $(y_1(c),y_2(c))$ and $(\tilde y_1(c),\tilde y_2(c))$ are linearly 
	dependent. Uniqueness of solutions gives that the functions $(y_1(a),y_2(a))$ and $(\tilde y_1(a),\tilde y_2(a))$ 
	are linearly dependent. Again evaluating the limit at $l_-$ yields that they are equal. 
\end{Remark}

\noindent
Another corollary of \Cref{M53} is that positivity is inherited.

\begin{corollary}
\label{M110}
	Let $p>-\frac 12$ and $(P,\psi)\in\bb R^{2\times 2}\times\bb R$. If $P\geq 0$, then 
	$\widehat\Xi_p(P,\psi)\in\HB$. 
\end{corollary}
\begin{proof}
	We actually are going to show that for all $a>0$ the function $a\odot_p A-ia\odot_p B$ is a Hermite-Biehler function. 
	Clearly, this is equivalent to the statement that all functions $E_a(z)\DE A(a,z)-iB(a,z)$ where $A(a,z)$ and $B(a,z)$ 
	are as in \cref{M65} are Hermite-Biehler functions. We know from \Cref{M104} that $E_a$ has no real zeroes, and we must
	prove positivity of the reproducing kernel $K_{E_a}$. 

	For $b>0$ let $W_b(a,z)$ be the unique solution of the initial value problem
	\begin{equation}\label{eq:2}
		\left\{
		\begin{array}{l}
			\frac{\partial}{\partial a}W_b(a,z)J=zW_b(a,z)H_{P,\psi}(a)\quad\text{for }a\in(0,\infty)
			,
			\\[2mm]
			W_b(b,z)=I
			.
		\end{array}
		\right.
	\end{equation}
	Then, by uniqueness of solutions, 
	\[
		\forall a\in(0,\infty)\DP
		\big(A(a,z),B(a,z)\big)=\big(A(b,z),B(b,z)\big)\cdot W_b(a,z)
		.
	\]
	From \eqref{eq:1} we obtain the kernel relation 
	\begin{align}
		K_{E_a}&(z,w)-K_{E_b}(z,w)=\nonumber\\
		&\big(A(b,z),B(b,z)\big)\frac{W_b(a,z)JW_b(a,w)^*-J}{z-\ov w}\big(A(b,w),B(b,w)\big)^*
		.\label{eq:3}
	\end{align}
	By \eqref{eq:2}, it follows for $b\leq a$ that
	\[
	W_b(a,z)JW_b(a,w)^*-J=(z-\ov w)\int_b^aW_b(c,z)H_{P,\psi}(c)W_b(c,w)^*dc.
	\]
	Our assumption that $P\geq 0$ implies that $H_{P,\psi}(a)\geq 0$ for all $a\in(0,\infty)$, and hence the kernel in \eqref{eq:3} is positive semidefinite for all $b\leq a$. By \Cref{M68} we have $\lim_{b\downarrow 0}K_{E_b}(z,w)=0$, and it follows that $K_{E_a}(z,w)$ is positive semidefinite. 
\end{proof}

%\newpage
%
%
%
\section[\textcolor{ForestGreen}{Solution of the canonical system via special functions}]{Solution of the canonical system 
	via special functions}
%ForestGreen / Dandelion / BrickRed / Fuchsia
%
%

\noindent
It is stated in \cite[p.205]{debranges:1962b} that ``homogeneous spaces of entire functions are related 
to Bessel functions and more general confluent hypergeometric functions'', however, ``it becomes tedious and awkward in handling
entire functions'', and thus actual formulae are not proven. 
In this section we provide the formulae for $\Xi_p(P,\psi)$ including all necessary computations, and things turn out not as 
awkward as one might expect.

Let us recall the definition of confluent hypergeometric (limit) functions:
\[
	M(\alpha,\beta,z)\DE\sum_{n=0}^\infty\frac{(\alpha)_n}{(\beta)_n}\cdot\frac{z^n}{n!}
	,\quad
	\prescript{}{0}{F}_1(\beta,z)\DE\sum_{n=0}^\infty\frac 1{(\beta)_n}\cdot\frac{z^n}{n!}
	,
\]
where $\alpha,z\in\bb C$ and $\beta\in\bb C\setminus(-\bb N_0)$. The symbol $(\Dummy)_n$ denotes the rising factorial, i.e., 
\[
	(\alpha)_0=1,\quad (\alpha)_{n+1}=(\alpha)_n(\alpha+n)\quad\text{for }n\in\bb N_0
	.
\]
The function $ \prescript{}{0}{F}_1$ is indeed a limit of $M$, namely, it holds that 
\begin{equation}
\label{M120}
	\prescript{}{0}{F}_1(\beta,z)=\lim_{|\alpha|\to\infty}M(\alpha,\beta,\smfrac z\alpha)
\end{equation}
locally uniformly in $z$ and $\beta$. 

By \Cref{M53} the solution of the canonical system with Hamiltonian $H_{P,\psi}$ is known once
$(A,B)=\Xi_p(P,\psi)$ has been computed. Our aim in this section is to prove the following explicit formulae. These formulae can be
derived by rewriting the differential equation corresponding to \cref{M54} as a scalar second-order differential equation of Kummer
type. A similar computation was carried out in \cite{langer.pruckner.woracek:asysupp-arXiv}. We shall not pursue this approach here.
Instead, we show that the coefficients of the power series expansion satisfy the required recurrence relation.  

\begin{Theorem}
\label{M36}
	Let $p\in\bb R\setminus(-\frac 12\bb N)$, 
	let $P=\smmatrix{\kappa_1}{\kappa_3}{\kappa_3}{\kappa_2}\in\bb R^{2\times 2}$ be a symmetric matrix, 
	and let $\psi\in\bb R$. 
	Let $\kappa\in\bb C$ be a square root of $\det P$, and set 
	\[
		\sigma\DE 2p\kappa_3-\psi\kappa_1,\qquad 
		\alpha\DE \smfrac\sigma{2i\kappa}+p\text{ if }\kappa\neq 0
		.
	\]
	As usual, we write $\Xi_p(P,\psi)=(A,B)$. 
	\begin{Enumerate}
	\item If $\det P\neq 0$, then we have 
		\begin{align}
		\nonumber
			A(z)= &\, 
			e^{i\kappa z}\Big[\frac 12 M(\alpha,2p+1,-2i\kappa z)+\frac 12 M(\alpha+1,2p+1,-2i\kappa z)
			\\
		\label{M123}
			&\, 
			-\frac{\kappa_3}{2p+1}zM(\alpha+1,2p+2,-2i\kappa z)\Big]
			,
			\\
		\label{M124}
			B(z)= &\,
			e^{i\kappa z}\frac{\kappa_1}{2p+1}zM(\alpha+1,2p+2,-2i\kappa z\big)
			.
		\end{align}
	\item If $\det P=0$, then we have 
		\begin{align}
		\label{M125}
			A(z)= &\, 
			\prescript{}{0}{F}_1(2p+1,-\sigma z)-\frac{\kappa_3}{2p+1}z\prescript{}{0}{F}_1(2p+2,-\sigma z)
			,
			\\
		\label{M126}
			B(z)= &\, 
			\frac{\kappa_1}{2p+1}z\prescript{}{0}{F}_1(2p+2,-\sigma z)
			.
		\end{align}
	\end{Enumerate}
\end{Theorem}

\noindent
Let us point out that the functions written in \Cref{M36} on the right sides of \cref{M123}--\cref{M126} do not depend on the 
choice of the square root $\kappa$. This is easy to check using the Kummer transformation \cite[\phantom{}13.2.39]{nist:2010}.

\begin{Remark}
\label{M127}
	The function $\prescript{}{0}{F}_1$ can be expressed in terms of Bessel functions of the first kind 
	\[
		J_\nu(z)\DE\sum_{n=0}^\infty\frac{(-1)^n}{n!\Gamma(n+\nu+1)}\Big(\frac z2\Big)^{2n+\nu}
		.
	\]
	The formula establishing this reads as
	\[
		J_\nu(z)=\frac{(\frac z2)^\nu}{\Gamma(\nu+1)}\prescript{}{0}{F}_1(\nu+1,-\smfrac{z^2}4)
		.
	\]
	Based on this fact, the formulae in the boundary case \Cref{M36}\,$(ii)$ could also be written in terms of Bessel
	functions.

	For a certain particular case, namely when $\beta=2\alpha$, the Kummer function $M$ is related to modified Bessel 
	functions 
	\[
		I_\nu(z)\DE\sum_{n=0}^\infty\frac{1}{n!\Gamma(n+\nu+1)}\Big(\frac z2\Big)^{2n+\nu}
		.
	\]
	The formula is
	\[
		I_\nu(z)=\frac{(\frac z2)^\nu}{\Gamma(\nu+1)}e^{-z}M(\nu+\smfrac 12,2\nu+1,2z)
		.
	\]
	If $\sigma=0$ this allows to rewrite the function from \cref{M124} to an expression involving only 
	modified Bessel functions. If $\sigma=0$ and $p\neq 0$, the same holds for the function from \cref{M123}. 
\end{Remark}
\begin{proof}[Proof of \Cref{M36}]
	
	 Proof of $(i)$:
	 
	 Let us set
	\begin{align*}
		\nonumber
		\tilde A(z)= &\, 
		\Big[\frac 12 M(\alpha,2p+1,-2i\kappa z)+\frac 12 M(\alpha+1,2p+1,-2i\kappa z)
		\\
		&\, 
		-\frac{\kappa_3}{2p+1}zM(\alpha+1,2p+2,-2i\kappa z)\Big]
		,
		\\
		\tilde B(z)= &\,
		\frac{\kappa_1}{2p+1}zM(\alpha+1,2p+2,-2i\kappa z\big)
		.
	\end{align*}
	and 
	\begin{equation*}
		\big(\tilde A(a,z),\tilde B(a,z)\big)\DE\big([a\odot_p \tilde A](z),[a\odot_p \tilde B](z)\big)\ms D_\psi(a)^{-1}
	\end{equation*}
	It is readily seen that 
	\begin{equation*}
		\big(\tilde A(a,z),\tilde B(a,z)\big)=e^{ia\kappa z}\big([a\odot_p A](z),[a\odot_p  B](z)\big)\ms D_\psi(a)^{-1}
	\end{equation*}
	
	Since $X(t)$ satisfies $X'(t)J=zX(t)H(t)$ if and only if $Y(t)=e^{i\kappa tz}X(t)$ satisfies $Y'(t)J=z(i\kappa J+H(t))$ we see that $\big(A(a,z),B(a,z)\big)$ satisfies \cref{M54}
	if and only if
	\begin{equation}
		\big(\tilde A(b,z),\tilde B(b,z)\big)J-\big(\tilde A(a,z),\tilde B(a,z)\big)J=z\int_a^b\big(A(c,z),B(c,z)\big)(-i\kappa J+ H_{P,\psi}(c))dc
		\label{M92}
	\end{equation}
	Setting $\tilde A(z)=\sum_{n=0}^\infty\tilde\alpha_nz^n$, $\tilde B(z)=\sum_{n=0}^\infty\tilde\beta_nz^n$ we get that $\tilde \alpha_0=1,\tilde \beta_0=0$ and for $n\geq 1$ we have
	\begin{align}\label{M93}
		\tilde \alpha_n= &\,-\frac{(-2i\kappa)^{n-1}(\alpha+1)_{(n-1)}(n\kappa_3+i\kappa(2\alpha+n))}{(2p+1)_nn!}
		\\
		\tilde\beta_n= &\,\frac{(-2i\kappa)^{n-1}\kappa_1(\alpha+1)_{n-1}}{(2p+1)_n(n-1)!}\label{M94}
		.
	\end{align}
	As in the proof of \Cref{M53}, $\big(\tilde A(a,z),\tilde B(a,z)\big)$ satisfies \cref{M92} if and only if $((\tilde\alpha_n,\tilde\beta_n))_{n\in\bb N_0}$ satisfy
	\begin{equation*}
		\left\{
		\begin{array}{l}
			(\tilde\alpha_{n+1},\tilde\beta_{n+1})=(\tilde\alpha_n,\tilde\beta_n)\cdot\frac{-1}{(n+1)(2p+n+1)}(PJ+i\kappa)\smmatrix{2p+n+1}{0}{\psi}{n+1}
			\quad\text{for }n\in\bb N_0
			,
			\\[4mm]
			(\tilde\alpha_0,\tilde\beta_0)=(1,0)
			.
		\end{array}
		\right.
	\end{equation*}
	
	Expanding the matrices, we see that the recursion is given by
	\begin{equation*}
		\begin{array}{l}
			(\tilde\alpha_{n+1},\tilde\beta_{n+1})=(\tilde\alpha_n,\tilde\beta_n)\cdot\frac{-1}{(n+1)(2p+n+1)}\smmatrix{(\kappa_3+i\kappa)(2p+n+1)-\kappa_1\psi}{-\kappa_1(n+1)}{\kappa_2(2p+n+1)+(-\kappa_3+i\kappa)\psi}{(-\kappa_3+i\kappa)(n+1)}.
		\end{array}
	\end{equation*}
	
	Let us start with the recursion for $\tilde\beta_{n+1}$: Plugging in \cref{M93} and \cref{M94} and using $(2p+n+1)(2p+1)_n=(2p+1)_{n+1}$ and $(\alpha+1)_{n-1}(\alpha+n)=(\alpha+1)_{n}$ we get
	\begin{align*}
		-\frac{1}{(2p+n+1)}&\left(-\tilde\alpha_n\kappa_1+\tilde\beta_n(-\kappa_3+i\kappa)\right)\\
		&=\frac{\kappa_1(-2i\kappa)^{n-1}(\alpha+1)_{n-1}}{(2p+1)_{n+1}(n-1)!}\left(\frac{n\kappa_3+i\kappa(2\alpha+n)}{n}+(-\kappa_3+i\kappa)\right)\\
		&=\frac{\kappa_1(-2i\kappa)^{n-1}(\alpha+1)_{n-1}}{(2p+1)_{n+1}(n-1)!}\left(\frac{(-i\kappa)(2\alpha+2n)}{n}\right)=\tilde\beta_{n+1}
	\end{align*}

	Let us now consider the recursion for $\tilde\alpha_{n+1}$: Using the recursion for $\tilde\beta_{n+1}$, the recursion for $\tilde\alpha_{n+1}$ simplifies to 
	\[
	-\frac{1}{n+1}(\tilde\alpha_n(\kappa_3+i\kappa)+\tilde\beta_n\kappa_2-\tilde\beta_{n+1}\psi)
	\]
	Straightforward manipulations as above using that $\kappa^2=\kappa_1\kappa_2-\kappa_3^2$ show that this is equal to 
	\[
	-\frac{(-2i\kappa)^n(\alpha+1)_n}{(2p+1)_{n+1}(n+1)!}((i\kappa+\kappa_3)(2p+n+1)-\psi\kappa_1).
	\]
	Thus, it remains to show that 
	\[
	(i\kappa+\kappa_3)(2p+n+1)-\psi\kappa_1=(n+1)\kappa_3+i\kappa(2\alpha+n+1).
	\]
	This can be easily checked from the definition of $\alpha$.

	Proof of $(ii)$: 	Assume first that $\sigma\neq 0$. The formulae \cref{M125}, \cref{M126} are obtained by passing to the limit 
	``$\kappa\to 0$'' in \cref{M123}, \cref{M124} using the formula \cref{M120}. 
	Note here that, clearly, the set of all invertible symmetric matrices is dense in the set of all symmetric matrices. 
	
	We have 
	\[
	\lim_{\kappa\to 0}(-2i\kappa)\alpha=\lim_{\kappa\to 0}(-2i\kappa)(\alpha+1)=-\sigma
	.
	\]
	Since $\sigma\neq 0$, we have $|\alpha|\to\infty$ when $\kappa\to 0$. Hence 
	\[
	\lim_{\kappa\to 0}M(\alpha,2p+1,-2i\kappa z)
	=\lim_{\kappa\to 0}M\big(\alpha,2p+1,\smfrac{(-2i\kappa)\alpha}\alpha\big)
	=\prescript{}{0}{F}_1(2p+1,-\sigma z)
	,
	\]
	and analogously
	\[
	\lim_{\kappa\to 0}M(\alpha+1,2p+1,-2i\kappa z)=\prescript{}{0}{F}_1(2p+1,-\sigma z)
	,
	\]
	\[
	\lim_{\kappa\to 0}M(\alpha+1,2p+2,-2i\kappa z)=\prescript{}{0}{F}_1(2p+2,-\sigma z)
	.
	\]
	This gives \cref{M125} and \cref{M126}. 
	
	The left and right sides of \cref{M125} and \cref{M126} depend continuously on $\sigma$, and hence the equality 
	holds also for $\sigma=0$. 
\end{proof}

\noindent
As a corollary of the explicit formulae in \Cref{M36} we obtain that $A,B$ are of bounded type in the upper and the lower half-plane:
\begin{corollary}
\label{M51}
	Let $p\in\bb R\setminus(-\frac 12\bb N)$, 
	let $P=\smmatrix{\kappa_1}{\kappa_3}{\kappa_3}{\kappa_2}\in\bb R^{2\times 2}$ be a symmetric matrix, 
	and let $\psi\in\bb R$. Assume that $\det P\geq 0$, and write as usual $\Xi_p(P,\psi)=(A,B)$. 
	Then $A$ and $B$ are of bounded type in the upper and the lower half-plane. 
\end{corollary}
\begin{proof}
	If $\det P=0$, we know that $A$ and $B$ are of order $\frac 12$, since the Bessel functions are of exponential type. 
	Assume that $\det P\neq 0$. 
	We know that $A$ and $B$ are of finite exponential type, cf. \cref{M35}. By Krein's Theorem (e.g.\ 
	\cite[Theorem~6.17]{rosenblum.rovnyak:1994}), it is thus enough to check
	convergence of the logarithmic integrals 
	\[
		\int_{-\infty}^\infty\frac{\log^+|A(x)|}{1+x^2}\DD x
		\quad\text{and}\quad
		\int_{-\infty}^\infty\frac{\log^+|B(x)|}{1+x^2}\DD x
		.
	\]
	We use the known asymptotics for confluent hypergeometric functions, see e.g.\ 
	\cite[\phantom{}13.5.1]{abramowitz.stegun:1964}:
	\[
		\frac{M(\alpha,\beta,x)}{\Gamma(\beta)}=
		\frac{e^{\pm i\pi\alpha}x^{-\alpha}}{\Gamma(\beta-\alpha)}\cdot\big(1+\Smallo({\textstyle\frac 1{|x|}})\big)
		+\frac{e^x x^{\alpha-\beta}}{\Gamma(\alpha)}\cdot\big(1+\Smallo({\textstyle\frac 1{|x|}})\big)
		.
	\]
	Since $\det P\geq 0$ we have $\kappa\in\bb R$. The above asymptotic expansions thus show that $|A(x)|$ and $|B(x)|$ are 
	bounded by some power for $x\in\bb R$.
\end{proof}

%\newpage
%
%
%
\section[\textcolor{ForestGreen}{Isometric inclusions of rescaled spaces}]{Isometric inclusions of rescaled spaces}
%ForestGreen / Dandelion / BrickRed / Fuchsia
%
%

It is a structural property of a space $\mc H(E)$ whether or not there are spaces $\mc H(a\odot_p E)$ which belong to its 
chain of de~Branges subspaces. 

\begin{Definition}
\label{M5}
	Let $p\in\bb R$ and $E\in\HB$. Then we denote 
	\[
		\mc O_p(E)\DE\big\{a\in\bb R^+ \DS \mc H(a\odot_p E)\subseteq\mc H(E)\text{ isometrically}\big\}
		.
	\]
\end{Definition}

\noindent
Trivially, $1\in\mc O_p(E)$. By \Cref{M3} we have $a\in\mc O_p(E)$ if and only if the map $F\mapsto a\odot_{p+\frac 12}F$ maps 
$\mc H(E)$ isometrically into itself. In particular, the set $\mc O_p(E)$ depends only on the space $\mc H(E)$ and not on 
the particular Hermite-Biehler function generating it. 

Since we have an underlying continuous group action, it can be expected that the set $\mc O_p(E)$ has some 
structure. To investigate it, we start with a basic fact.

\begin{lemma}
\label{M7}
	Let $p\in\bb R$ and $E\in\HB$, and let further $a,b\in\bb R^+$. Then the following statements hold. 
	\begin{Enumerate}
	\item ${\displaystyle
		\begin{aligned}[t]
			& \mc H(a\odot_p E)\subseteq\mc H(b\odot_p E)\text{ isometrically}
			\\
			& \mkern50mu\Longleftrightarrow\quad 
			\forall c\in\bb R^+ \DP \mc H(ca\odot_p E)\subseteq\mc H(cb\odot_p E)\text{ isometrically}
			\\
			& \mkern50mu\Longleftrightarrow\quad 
			\frac ab\in\mc O_p(E)
			.
		\end{aligned}
		}$
	\item The set $\mc O_p(E)$ is a subsemigroup of $\bb R^+$. 
	\item ${\displaystyle
		\forall a\in\bb R^+\DP \mc O_p(a\odot_p E)=\mc O_p(E)
		.
		}$
	\end{Enumerate}
\end{lemma}
\begin{proof}
	For the proof of $(i)$ let $a,b,c>0$. 
	By \Cref{M3} the map $F\mapsto c\odot_{p+\frac 12} F$ is an isometric isomorphism of 
	$\mc H(a\odot_p E)$ onto $\mc H(c\odot_p (a\odot_p E))$. Now note that 
	$c\odot_p (a\odot_p E)=(ca)\odot_p E$. The same holds for $b$ in place of $a$, and we see that the first condition 
	implies the second. The converse implication is trivial; just use ``$c=1$''. 
	To prove that the first and third conditions are equivalent, apply the already proven with $c\DE\frac 1b$ for 
	``$\Rightarrow$'' and with $c\DE b$ for ``$\Leftarrow$''.

	To show that $\mc O_p(E)$ is closed under multiplication, let $a,b\in\mc O_p(E)$. 
	Then $\mc H(b\odot_p E)\subseteq\mc H(E)$ isometrically, and hence also 
	$\mc H(ab\odot_p E)\subseteq\mc H(a\odot_p E)$ isometrically. Since $\mc H(a\odot_p E)\subseteq\mc H(E)$ 
	isometrically, it follows in turn that $\mc H(ab\odot_p E)\subseteq\mc H(E)$ isometrically.

	Finally, to show $(iii)$, note that by the already proven equivalence in $(i)$ we have 
	\begin{align*}
		c\in\mc O_p(a\odot_p E)\quad\Longleftrightarrow\quad & 
		\mc H\big(\underbrace{c\odot_p(a\odot_p E)}_{=ca\odot_p E}\big)\subseteq\mc H(a\odot_p E)\text{ isometrically}
		\\
		\Longleftrightarrow\quad & c\in\mc O_p(E)
		,
	\end{align*}
	and the same for $b$ in place of $a$. 
\end{proof}

\noindent
Our next statement is that generically isometric inclusions can occur only for $p>-\frac 12$ and $a\leq 1$. 
The case of one-dimensional spaces is exceptional and corresponds to the boundary case $p=-\frac 12$.

\begin{proposition}
\label{M6}
	Let $E\in\HB$. Then the following statements hold.
	\begin{Enumerate}
	\item If $\dim\mc H(E)>1$ and $p\in\bb R$ with $\mc O_p(E)\neq\{1\}$, then $p>-\frac 12$ and 
		$\mc O_p(E)\subseteq(0,1]$. 
	\item If $\dim\mc H(E)=1$ and $p\in\bb R$ with $\mc O_p(E)\neq\{1\}$, then 
		$p=-\frac 12$ and $\mc H(E)=\Span\{1\}$. 
	\item If $\mc H(E)=\Span\{1\}$, then $\mc O_{-\frac 12}(E)=\bb R^+$.
	\end{Enumerate}
\end{proposition}
\begin{proof}
	For the proof of $(i)$ assume that $\dim\mc H(E)>1$. In the first step we prove that 
	\begin{equation}
	\label{M9}
		\forall p\in\bb R \DP \mc O_p(E)\subseteq(0,1]
	\end{equation}
	Choose a phase function $\varphi_E$ associated with $E$ (cf.\ \Cref{M81}). By 
	\cite[Theorem~22,Problem~46]{debranges:1968} our assumption that $\mc H(E)$ is at least two-dimensional implies that 
	\[
		\lim_{x\to\infty}\big(\varphi_E(x)-\varphi_E(-x)\big)>\pi
		.
	\]
	Hence, we find $x_0>0$ with $\varphi_E(x_0)-\varphi_E(-x_0)=\pi$. 

	Let $p\in\bb R$ and $a\in\mc O_p(E)$. Clearly, the function $x\mapsto\varphi_E(ax)$ is a phase function associated with
	$a\odot_p E$. Since $\mc H(a\odot_p E)\subseteq\mc H(E)$ isometrically, \cite[Problem~93]{debranges:1968} yields 
	\[
		\varphi_E(ax_0)-\varphi_E(-ax_0)=\varphi_{a\odot_p E}(x_0)-\varphi_{a\odot_p E}(-x_0)
		\leq\varphi_E(x_0)-\varphi_E(-x_0)
		.
	\]
	Since $\varphi_E$ is strictly increasing, it follows that $a\leq 1$. The proof of \cref{M9} is complete. 

	As a consequence, we obtain that 
	\begin{equation}
	\label{M10}
		\forall p\in\bb R \DP \mc O_p(E)\neq\{1\}\ \Rightarrow\ p\geq-\frac 12
	\end{equation}
	To see this, recall the kernel relation \cref{M4}. It implies in particular that 
	\[
		K_{a\odot_p E}(0,0)=a^{2p+1}K_E(0,0).
	\]
	If $\mc H(a\odot_p E)\subseteq\mc H(E)$ isometrically, the reproducing kernels satisfy 
	$K_{a\odot_p E}(0,0)\leq K_E(0,0)$, since these quantities are the (square of the) norm of point evaluation at $0$. 
	It follows that $a^{2p+1}\leq 1$. Knowing that $a$ cannot exceed $1$, \cref{M10} follows. 

	It remains to exclude the case that $p=-\frac 12$. 	Assume that $a\in\mc O_{-\frac 12}(E)$. By \Cref{M7} 
	we have the chain of isometric inclusions 
	\begin{equation}
	\label{M11}
		\mc H(a^2\odot_{-\frac 12} E)\subseteq\mc H(a\odot_{-\frac 12} E)\subseteq\mc H(E)
		.
	\end{equation}
	Since $K_{a^2\odot_{-\frac 12} E}(0,0)=K_E(0,0)$, i.e., the interval between these two de~Branges subspaces is
	indivisible of type $\frac\pi2$, we have
	\[
		\dim\Big(\raisebox{2pt}{$\mc H(E)$}\Big/\raisebox{-3pt}{$\mc H(a^2\odot_{-\frac 12} E)$}\Big)\leq 1
		.
	\]
	In the above chain \cref{M11} thus at least one inclusion must hold with equality. Using the appropriate isometry 
	($F\mapsto\frac 1a\odot_0F$ or $F\mapsto a\odot_0F$), this implies that equality holds
	throughout. Hence, certainly $\mc H(a\odot_{-\frac 12} E)=\mc H(E)$, and we see that 
	$F\mapsto a\odot_0F$ is an isometric bijection of $\mc H(E)$ onto itself. Its inverse, which is 
	$F\mapsto \frac 1a\odot_0F$, thus has the same property. Now \cref{M9} implies that $a=1$. 
	The proof of $(i)$ is complete. 

	We come to the proof of $(ii)$. Assume that $\mc H(E)=\Span\{G\}$ with some entire function $G$, and assume further that
	$p\in\bb R$ and $a\in\mc O_p(E)\setminus\{1\}$. The function $a\odot_{p+\frac 12}G$ belongs to the space $\mc H(E)$ and
	has the same norm as $G$. Thus, $a\odot_{p+\frac 12}G=G$. Writing the power series expansion of $G$ as 
	$G(z)=\sum_{n=0}^\infty\gamma_n z^n$, this gives 
	\[
		\sum_{n=0}^\infty a^{p+\frac 12+n}\gamma_n z^n=\sum_{n=0}^\infty\gamma_n z^n
		,
	\]
	and comparing coefficients yields 
	\[
		\forall n\in\bb N_0 \DP \big(a^{p+\frac 12+n}-1\big)\gamma_n=0
	\]
	Since $a\neq 1$ and $\gamma_0=G(0)\neq 0$, recall here the fourth property in \Cref{M101}, this implies that 
	\[
		p=-\frac 12\quad\text{and}\quad\forall n\geq 1 \DP \gamma_n=0
		.
	\]
	Thus, $G$ is constant and $\mc H(E)=\Span\{1\}$. 

	For the proof of $(iii)$ it suffices to note that each map $F\mapsto a\odot_0 F$ acts as the identity on 
	$\Span\{1\}$. 
\end{proof}

\noindent
Now the structure of the set $\mc O_p(E)$ can be clarified. 

\begin{proposition}
\label{M8}
	Let $p>-\frac 12$ and $E\in\HB$. Then one of the following statements holds.
	\begin{Itemize}
	\item $\mc O_p(E)=\{1\}$;
	\item $\exists a_0\in(0,1)\DP \mc O_p(E)=\{(a_0)^n\DS n\in\bb N_0\}$;
	\item $\mc O_p(E)=(0,1]$. 
	\end{Itemize}
\end{proposition}
\begin{proof}
	We have already seen that $\mc O_p(E)$ is a subsemigroup of $(0,1]$. The first step in the proof is to show that 
	$\mc O_p(E)$ is also invariant under suitable quotients:
	\begin{equation}
	\label{M12}
		\forall a,b\in\mc O_p(E),a\leq b \DP \frac ab\in\mc O_p(E)
		.
	\end{equation}
	Let $a,b\in\mc O_p(E)$. If $a=b$, there is nothing to prove, hence assume that $a<b$. 
	Both spaces $\mc H(a\odot_p E)$ and $\mc H(b\odot_p E)$ are contained isometrically in $\mc H(E)$. By the 
	Ordering Theorem (e.g.\ \cite[Theorem~35]{debranges:1968}), either 
	$\mc H(a\odot_p E)\subseteq\mc H(b\odot_p E)$ or $\mc H(b\odot_p E)\subseteq\mc H(a\odot_p E)$. Since 
	\[
		K_{a\odot_p E}(0,0)=a^{2p+1}K_E(0,0)<b^{2p+1}K_E(0,0)=K_{b\odot_p E}(0,0)
		.
	\]
	the first case must take place. Now \Cref{M7} applies, and we obtain that $\frac ab\in\mc O_p(E)$. 

	Knowing \cref{M12}, it follows that the set $\mc G\DE\mc O_p(E)\cup\mc O_p(E)^{-1}$ is a subgroup of $\bb R^+$ 
	and satisfies $\mc G\cap(0,1]=\mc O_p(E)$. A subgroup of $\bb R^+$ is either 
	\begin{Itemize}
	\item trivial, i.e., $\mc G=\{1\}$, or 
	\item nontrivial and cyclic, i.e., $\mc G=\{(a_0)^n\DS n\in\bb Z\}$ for some $a_0\in(0,1)$, or 
	\item dense in $\bb R^+$. 
	\end{Itemize}
	To complete the proof of the present assertion, it is enough to note that $\mc O_p(E)$ is closed in $(0,1]$:
	By the formula \cref{M4} for $K_{a\odot_p E}$ we have $\lim_{n\to\infty}K_{a_n\odot_p E}=K_{a\odot_p E}$ locally
	uniformly. Now remember \Cref{M86} and observe that $K_{a\odot_p E}(0,0)\neq 0$.
\end{proof}

\noindent
Let us show by examples that each of the three alternatives can occur. 

\begin{example}
\label{M69}
	Let $\mc H(E)$ be a de~Branges space with $\dim\mc H(E)>1$ and $1\in\mc H(E)$. We assert that $\mc O_p(E)=\{1\}$ for all
	$p$. To see this, assume towards a contradiction that there exists $a\in\mc O_p(E)\setminus\{1\}$. 
	We have $a\odot_{p+\frac 12}1=a^{p+\frac 12}$, and hence 
	\[
		\|1\|=\|a\odot_{p+\frac 12}1\|=a^{p+\frac 12}\cdot\|1\|
		.
	\]
	This contradicts the fact that $p>-\frac 12$ by \Cref{M6}.
\end{example}

\begin{example}
\label{M82}
	Let $\mr H$ be a Hamiltonian on $(0,1)$ with $\Tr\mr H=1$ such that the interval $(0,1)$ is not indivisible, 
	and define a Hamiltonian $H$ on $(0,\infty)$ by 
	\[
		H(t)\DE \mr H\Big(\frac t{2^n}-1\Big)
		\quad\text{for }n\in\bb Z,t\in(2^n,2^{n+1})
		.
	\]
	Then $H$ is multiplicatively periodic with period $2$, i.e., $H(2t)=H(t)$ for all $t>0$. 
	Let $(A(t,z),B(t,z))$ be the unique solution of the initial value problem
	\[
		\left\{
		\begin{array}{l}
			\frac{\partial}{\partial t}(A(t,z),B(t,z))J=z(A(t,z),B(t,z))H(t)\quad\text{for }t>0
			,
			\\[2mm]
			(A(0,z),B(0,z))=(1,0)
			.
		\end{array}
		\right.
	\]
	Since $H$ is periodic, uniqueness of the solution implies that 
	\[
		\forall t\geq 0,z\in\bb C\DP 
		\big(A(2t,z),B(2t,z)\big)=\big(A(t,2z),B(t,2z)\big)
		.
	\]
	In other words, the functions $E(t,z)\DE A(t,z)-iB(t,z)$ satisfy 
	\[
		\forall t\geq 0\DP
		E(2t,.)=2\odot_0 E(t,.)
		.
	\]
	Choose $t_0\in(0,1)$ such that $t_0$ is not inner point of an $\mr H$-indivisible interval. 
	Then each point $2^n(t_0+1)$ with $n\in\bb Z$ is not inner point of an $H$-indivisible interval, and hence 
	\[
		\forall n,m\in\bb Z,n\leq m\DP
		\mc H\big(E(2^n(t_0+1),.)\big)\subseteq\mc H\big(E(2^m(t_0+1),.)\big)\text{ isometrically}
		.
	\]
	Since $E(2^n(t_0+1),.)=2^n\odot_0 E(t_0+1,.)$, we see that 
	\[
		\{2^n\DS n\in -\bb N_0\}\subseteq\mc O_0(E(t_0+1,.))
		.
	\]
	Now we want to impose some condition on $\mr H$ which guarantees that 
	\[
		\mc O_0(E(t_0+1,.))\neq(0,1]
		.
	\]
	For example, assuming that $\mr H$ consists of a sequence of indivisible intervals does the job. 
	Assume that we have $N\geq 2$ and 
	\[
		0=t_0<t_1<\ldots<t_{N-1}<t_N=1
		,
	\]
	such that each interval $(t_{i-1},t_i)$ is $\mr H$-indivisible and that the types of each two successive intervals
	are different. Then all but countably many points of $(0,\infty)$ are inner point of some $H$-indivisible interval,
	and hence the chain of de~Branges subspaces of $\mc H(E(t_0+1,.))$ is countable. Because of \cref{M4} the map 
	$a\mapsto\mc H(a\odot_0 E(t_0+1,.))$ maps $\mc O_0(E(t_0+1,.))$ injectively into this chain.
	Thus $\mc O_0(E(t_0+1,.))$ cannot be equal to $(0,1]$. 
\end{example}

\begin{example}
\label{M13}
	The classical Paley-Wiener space $\PW_1$ is generated by the Hermite-Biehler function $E(z)\DE e^{-iz}$. 
	Let $F\in\mc H(E)$. For all $p\in\bb R$ and $a\in(0,1)$ the function $a\odot_{p+\frac 12}F$ belongs to $\mc H(E)$. 
	However, its norm coincides with the norm of $F$, if and only if $p=0$:
	\[
		\|a\odot_{p+\frac 12}F\|_{\mc H(E)}^2=\int_{\bb R}|a^{p+\frac 12}F(ax)|^2\DD x=
		\int_{\bb R}a^{2p}|F(y)|^2\DD y
		.
	\]
	Thus 
	\[
		\mc O_p(E)=
		\begin{cases}
			(0,1] \CAS p=0
			,
			\\[1mm]
			\{1\} \CASO
			.
		\end{cases}
	\]
\end{example}

\subsection*{Homogeneous de~Branges spaces}

Let us now turn to homogeneous spaces in the sense of \Cref{M70}. Note that a space $\mc H(E)$ is homogeneous of order 
$\nu$ (with some $\nu>-1$), if and only if $\mc O_{\nu+\frac 12}(E)=(0,1]$. 
We had decided to stick to the terminology introduced by L.de~Branges in \cite[Chapter~50]{debranges:1968} in \Cref{M70}, 
and this is why we have a shift by $\frac 12$ in the parameters, i.e., $\nu$ from that definition and $p$ from above 
are related as $p=\nu+\frac 12$.

For $E\in\HB$ such that $\mc H(E)$ is homogeneous of order $\nu$, the chain of de~Branges spaces 
\[
	\big\{\mc H(a\odot_{\nu+\frac 12}E) \DS a\in(0,1]\big\}
\]
is contained in the chain of all de~Branges subspaces of $\mc H(E)$. 
It follows from continuity that it exhausts this chain (cf.\ \Cref{M15} below). 
To prove this, we recall the following general and probably folklore fact.

\begin{lemma}
\label{M14}
	Let $\Omega$ be a nonempty set, $\mc I\subseteq\bb R$ an interval, and $(\mc H_t)_{t\in\mc I}$ a family of reproducing
	kernel Hilbert spaces of complex valued functions on $\Omega$. Denote the reproducing kernel of $\mc H_t$ as $K_t$. 
	Assume that 
	\begin{align*}
		& \forall s,t\in\mc I \DP t\leq s\ \Rightarrow\ \mc H_t\subseteq\mc H_s\text{ isometrically}
		,
		\\
		& \forall w\in\Omega \DP t\mapsto K_t(w,w)\text{ is continuous}
		.
	\end{align*}
	Then 
	\begin{align}
	\label{M28}
		& \forall s\in\mc I\setminus\{\sup\mc I\} \DP 
		\bigcap_{t\in\mc I,t>s}\mc H_t=\mc H_s
		,
		\\
	\label{M29}
		& \forall s\in\mc I\setminus\{\inf\mc I\} \DP 
		\Clos\bigcup_{t\in\mc I,t<s}\mc H_t=\mc H_s
		.
	\end{align}
\end{lemma}
\begin{proof}
	For $s,t\in\mc I$ with $s\leq t$ we denote by $P^t_s$ the orthogonal projection of $\mc H_t$ onto $\mc H_s$. 

	To prove \cref{M28}, let $s\in\mc I\setminus\{\sup\mc I\}$ be given. Choose $t_0\in\mc I$ with $t_0>s$, 
	denote $\mc H_{s+}\DE\bigcap_{t\in\mc I,t>s}\mc H_t$ and let $K_{s+}$ be the reproducing kernel of $\mc H_{s+}$. 
	The limit $\lim_{t\downarrow s}P^{t_0}_t$ exists in the strong operator topology and is the orthogonal projection of 
	$\mc H_{t_0}$ onto $\mc H_{s+}$. Thus 
	\[
		\forall w\in\Omega \DP 
		K_{s+}(.,w)=\lim_{t\downarrow s}K_t(.,w)
	\]
	in norm. We conclude that, for each $w\in\Omega$, 
	\[
		\|K_{s+}(.,w)\|=\lim_{t\downarrow s}\|K_t(.,w)\|=
		\lim_{t\downarrow s}K_t(w,w)^{\frac 12}=K_s(w,w)^{\frac 12}=\|K_s(.,w)\|
		.
	\]
	Since $K_s(.,w)$ is the orthogonal projection of $K_{s+}(.,w)$ onto $\mc H_s$, it follows that 
	$K_{s+}(.,w)=K_s(.,w)$. We see that $\mc H_{s+}=\mc H_s$. 

	The proof of \cref{M29} is carried out in the same way. 
\end{proof}

\begin{corollary}
\label{M15}
	Let $E\in\HB$, $\nu>-1$, and assume that $\mc H(E)$ is homogeneous of order $\nu$. Then 
	the chain of de~Branges subspaces of $\mc H(E)$ is equal to 
	\[
		\big\{\mc H(a\odot_{\nu+\frac 12}E) \DS a\in(0,1]\big\}
		.
	\]
\end{corollary}
\begin{proof}
	To shorten notation, set $p\DE\nu+\frac 12$. 
	Let $\mc H(\tilde E)$ be a de~Branges subspace of $\mc H(E)$. For each $a\in(0,1]$ we must have either 
	$\mc H(\tilde E)\subseteq\mc H(a\odot_p E)$ or $\mc H(a\odot_p E)\subseteq\mc H(\tilde E)$. 
	If $K_{a\odot_p E}(0,0)>K_{\tilde E}(0,0)$ or $a=1$, then the first case must take place. On the other hand, the second
	case must take place whenever $K_{a\odot_p E}(0,0)<K_{\tilde E}(0,0)$, and this inequality holds for all sufficiently
	small $a$ since $\lim_{a\downarrow 0}K_{a\odot_p E}(z,w)=0$ by \cref{M4}.
	Set 
	\[
		b\DE\Big(\frac{K_{\tilde E}(0,0)}{K_E(0,0)}\Big)^{\frac 1{2p+1}}\in(0,1]
		.
	\]
	The function $a\mapsto a\odot_p E$ is continuous, and \Cref{M14} implies 
	\begin{multline*}
		\mc H(b\odot_p E)=\Clos\bigcup_{a\in(0,b)}\mc H(a\odot_p E)\subseteq\mc H(\tilde E)
		\\
		\subseteq\mc H(E)\cap\bigcap_{a\in(b,1]}\mc H(a\odot_p E)=\mc H(b\odot_p E)
		.
	\end{multline*}
\end{proof}

\begin{corollary}
\label{M71}
	Let $\mc H$ be a de~Branges space and $\nu>-1$. If $\mc H$ is homogeneous of order $\nu$, then every de~Branges 
	subspace of $\mc H$ is homogeneous of order $\nu$. 
\end{corollary}
\begin{proof}
	Choose $E\in\HB$ such that $\mc H=\mc H(E)$. 
	We know that each de~Branges subspace of $\mc H(E)$ is of the form $\mc H(a\odot_{\nu+\frac 12}E)$ with some $a\leq 0$. 
	Now remember \Cref{M7}\,$(iii)$.
\end{proof}

%\newpage
%
%
%
\section[\textcolor{ForestGreen}{Structure of homogeneous spaces}]{Structure of homogeneous spaces}
%ForestGreen / Dandelion / BrickRed / Fuchsia
%
%

In this section we give the connection between homogeneous de~Branges spaces, canonical systems with Hamiltonians of the form 
$H_{P,\psi}$, and recurrence relations of the form \cref{M30}.

\begin{Definition}
\label{M32}
	Let $p>-\frac 12$. Then we set 
	\[
		\bb P_p\DE 
		\Big\{(P,\psi)\in\bb R^{2\times 2}\times\bb R \DS 
		P\geq 0,\ \binom 10,\binom{-\psi}{2p}\notin\ker P\Big\}
		\quad\text{if }p\neq 0
		,
	\]
	\begin{align*}
		\bb P_0\DE &\, 
		\Big\{(P,\psi)\in\bb R^{2\times 2}\times\bb R \DS 
		P\geq 0,\ \ker P=\{0\},\ \psi=0\Big\}
		\\
		&\, \cup
		\Big\{(P,\psi)\in\bb R^{2\times 2}\times\bb R \DS 
		P\geq 0,\ \binom 10\notin\ker P,\ \psi\neq 0\Big\}
		.
	\end{align*}
\end{Definition}

\noindent
The conditions occurring in this definition stem from \Cref{M41}. Also note that $\bb P_p\subseteq\bb P$ for all 
$p>-\frac 12$. 

\begin{Theorem}
\label{M16}
	Let $\nu>-1$ and set $p\DE\nu+\frac 12$. Then the following statements hold.
	\begin{Enumerate}
	\item Let $(P,\psi)\in\bb P_p$. Then 
		\begin{Ilist}
		\item $\widehat\Xi_p(P,\psi)$ is a Hermite-Biehler function with value $1$ at the origin, 
		\item the structure Hamiltonian of $\widehat\Xi_p(P,\psi)$ is a reparameterisation of $H_{P,\psi}|_{(0,1]}$ 
			(prolongued with an indivisible interval of type $\frac\pi 2$ up to $-\infty$ if 
			$-\frac 12<p<\frac 12$), 
		\item $\mc H(\widehat\Xi_p(P,\psi))$ is homogeneous of order $\nu$. 
		\end{Ilist}
	\item Let $E\in\HB$ with $E(0)=1$ be such that $\mc H(E)$ is homogeneous of order $\nu$. 
		Then there exists $(P,\psi)\in\bb P_p$, such that $E=\widehat\Xi_p(P,\psi)$ and that the structure 
		Hamiltonian $H_E$ is a reparameterisation of $H_{P,\psi}|_{(0,1]}$ 
		(prolongued by an indivisible interval if necessary). 
	\item Let $(P,\psi),(\tilde P,\tilde\psi)\in\bb P_p$, and write $P=\smmatrix{\kappa_1}{\kappa_3}{\kappa_3}{\kappa_2}$
		and $\tilde P=\smmatrix{\tilde\kappa_1}{\tilde\kappa_3}{\tilde\kappa_3}{\tilde\kappa_2}$. Then 
		$\mc H(\widehat\Xi_p(P,\psi))=\mc H(\widehat\Xi_p(\tilde P,\tilde\psi))$ isometrically, if and only if
		\begin{Ilist}
		\item $\kappa_1=\tilde\kappa_1$, 
		\item $\det P=\det\tilde P$, 
		\item $\psi-\tilde\psi=
			\frac{2p}{\kappa_1}\big[\kappa_3-\tilde\kappa_3\big]$.
		\end{Ilist}
	\end{Enumerate}
\end{Theorem}

\noindent
The proofs of items $(i)$ and $(iii)$ of \Cref{M16} are rather easy: the first is obtained by plugging together what we have 
already shown, and the latter by calculation. The essential part of the proof is to establish the existence result stated in 
$(ii)$.

\begin{proof}[Proof of \Cref{M16}\,$(i)$]
	We know from \Cref{M110} that $\widehat\Xi_p(P,\psi)\in\HB$. 
	Let us observe some properties of $H_{P,\psi}$. 
	It is clear that $H_{P,\psi}(a)\geq 0$ for all $a>0$ and that $H_{P,\psi}$ is locally integrable. The definition of 
	the class $\bb P_p$ and \Cref{M60} imply that $H_{P,\psi}\notin L^1((1,\infty),\bb R^{2\times 2})$ and that 
	$H_{P,\psi}\in L^1((0,1),\bb R^{2\times 2})$ if and only if $p\in(-\frac 12,\frac 12)$. Further $H_{P,\psi}$ contains no
	indivisible intervals: if $\Rank P=2$ this is trivial, and if $\Rank P=1$ it follows from \Cref{M41} by the definition 
	of $\bb P_p$.

	Let $A(a,z)$ and $B(a,z)$ be as in \cref{M65}. 
	\Cref{M53} implies that for all $a\in(0,1]$ the space $\mc H(A(a,z)-iB(a,z))$ is isometrically contained in 
	$\mc H(\widehat\Xi_p(P,\psi))$. Note here that changing the functions $A,B$ with a matrix from $\SL(2,\bb R)$ does not
	change the generated de~Branges space. However, by the definition \cref{M65} we have 
	$\mc H(A(a,z)-iB(a,z))=\mc H(a\odot_p\widehat\Xi_p(P,\psi))$ isometrically. Moreover, the function 
	\[
		\binom 10^*H_{P,\psi}(a)\binom 10=a^{2p}\binom 10^*P\binom 10
		,
	\]
	is integrable at $0$. Altogether we see that the structure Hamiltonian of $\widehat\Xi_p(P,\psi)$ is a
	reparameterisation of $H_{P,\psi}|_{(0,1]}$ (extended by an indivisible interval if necessary), and 
	that $\mc H(\widehat\Xi_p(P,\psi))$ is homogeneous of order $\nu$. 
\end{proof}

\begin{proof}[Proof of \Cref{M16}\,$(ii)$]
	Let $\mc H$ be a de~Branges space which is homogeneous of order $\nu$, and choose $E\in\HB$ with $E(0)=1$ such that 
	$\mc H=\mc H(E)$. 
	As usual, let $H_E$ denote the structure Hamiltonian of $E$, let $W_E(t,s,z)$ for $-\infty<t\leq s\leq 0$ be the 
	corresponding family of transfer matrices, and $E(t,z)$ for $t\in(-\infty,0]$ be the corresponding 
	family of Hermite-Biehler functions. Moreover, set 
	\[
		t_-\DE\sup\big\{t\in(-\infty,0] \DS (-\infty,t)\text{ $H_E$-indivisible of type }\frac\pi 2\big\}
		.
	\]
	\begin{Elist}
	\item By \Cref{M15} the chain of de~Branges subspaces of $\mc H(E)$ is $\{\mc H(a\odot_p E)\DS a\in(0,1]\}$, and by 
		\Cref{M14} this chain has no one-dimensional gaps. Hence, the Hamiltonian $H_E|_{(t_-,0)}$ contains no 
		indivisible intervals. In particular, the function 
		\begin{equation}
		\label{M18}
			\alpha(t)\DE\bigg(
			\frac 1{K_E(0,0)}\int\limits_{-\infty}^t\binom 10^*H_E(t)\binom 10 \DD t
			\bigg)^{\frac 1{2p+1}}
		\end{equation}
		is an absolutely continuous increasing bijection of $[t_-,0]$ onto $[0,1]$. Its inverse function is thus a 
		continuous increasing bijection of $[0,1]$ onto $[t_-,0]$. Let us point out that at the present stage we do not 
		know that $\alpha^{-1}$ is absolutely continuous. 

		Remembering \cref{M4}, we see that 
		\[
			\forall t\in(t_-,0] \DP \mc H(E(t,.))=\mc H(\alpha(t)\odot_p E)
			.
		\]
		Now define a function $\Psi\DF(0,1]\to\SL(2,\bb R)$ by letting $\Psi(a)$ be the unique matrix with
		\begin{equation}
		\label{M66}
			(a\odot_p A,a\odot_p B)=\big(A(\alpha^{-1}(a),.),B(\alpha^{-1}(a),.)\big)\Psi(a)
			,
		\end{equation}
		cf.\ \Cref{M99}.
		Clearly, $\Psi(1)=I$. Evaluating \cref{M66} and the relation obtained by differentiating \cref{M66} w.r.t.\ $z$ 
		at the point $z=0$, yields that $\Psi(a)$ is explicitly given as (here and in the following a prime denotes 
		differentiation w.r.t.\ $z$)
		\[
			\Psi(a)=
			\begin{pmatrix}
				1 & 0
				\\
				A'(\alpha^{-1}(a),0) & B'(\alpha^{-1}(a),0)
			\end{pmatrix}
			^{-1}
			\begin{pmatrix}
				a^p & 0
				\\
				a^{p+1}A'(0) & a^{p+1}B'(0)
			\end{pmatrix}
			.
		\]
		This formula together with the fact that $\Psi(a)\in\SL(2,\bb R)$ shows that $\Psi$ is continuous and triangular 
		of the form 
		\begin{equation}
		\label{M21}
			\Psi(a)=
			\begin{pmatrix}
				a^p & 0
				\\
				\text{\ding{100}} & a^{-p}
			\end{pmatrix}
			.
		\end{equation}
	\item In this step we consider the function
		\[
			W(a,b,z)\DE W_E(\alpha^{-1}(a),\alpha^{-1}(b),z)\quad\text{for }0<a\leq b\leq 1
			,
		\]
		and show the central relation 
		\begin{multline}
			\label{M19}
			\forall 0<a\leq b\leq 1,0<c\leq\frac 1b \DP 
			\\
			\Psi(a)^{-1}W(a,b,cz)\Psi(b)=\Psi(ca)^{-1}W(ca,cb,z)\Psi(cb)
			.
		\end{multline}
		Let $0<a\leq b\leq 1$. Then
		\begin{align}
		\label{M33}
			\big([a\odot_p A &\, ](z),[a\odot_p B](z)\big)\cdot\Psi(a)^{-1}W(a,b,z)\Psi(b)
			\\
		\nonumber
			= &\, \big(A(\alpha^{-1}(a),z),B(\alpha^{-1}(a),z)\big) W_E(\alpha^{-1}(a),\alpha^{-1}(b),z)\Psi(b)
			\\
		\nonumber
			= &\, \big(A(\alpha^{-1}(b),z),B(\alpha^{-1}(b),z)\big)\Psi(b)
			=\big([b\odot_p A](z),[b\odot_p B](z)\big)
			.
		\end{align}
		Let additionally $0<c\leq\frac 1b$, then we can compute 
		\begin{align*}
			\big([ca\odot_p A &\, ](z),[ca\odot_p B](z)\big)\cdot\Psi(a)^{-1}W(a,b,cz)\Psi(b)
			\\
			= &\, c\odot_p\Big[\big([a\odot_p A](z),[a\odot_p B](z)\big)\Psi(a)^{-1}W(a,b,z)\Psi(b)\Big]
			\\
			= &\, c\odot_p\Big[\big([b\odot_p A](z),[b\odot_p B](z)\big)\Big]
			=\big([cb\odot_p A](z),[cb\odot_p B](z)\big)
			\\
			= &\, \big([ca\odot_p A](z),[ca\odot_p B](z)\big)\cdot\Psi(ca)^{-1}W(ca,cb,z)\Psi(cb)
			.
		\end{align*}
		Uniqueness of the transfer matrix \cite[Problem~100]{debranges:1968} implies \cref{M19}. 
	\item We exploit \cref{M19} to determine $\Psi$. Using \cref{M19} with $b=1$ and evaluating at $z=0$ leads to 
		\[
			\forall a,c\in(0,1] \DP \Psi(ca)=\Psi(c)\Psi(a)
			.
		\]
		The function 
		\[
			\left\{
			\begin{array}{rcl}
				[0,\infty) & \to & \bb C^{2\times 2}
				\\
				x & \mapsto & \Psi(e^{-x})
			\end{array}
			\right.
		\]
		is a continuous semigroup of matrices, and hence is represented as the exponential of its infinitesimal 
		generator:
		\[
			\Psi(e^{-x})=\exp\big(Gx\big)\quad\text{where}\quad
			G\DE\lim_{x\downarrow 0}\frac{\Psi(e^{-x})-I}x
			.
		\]
		Since $\Psi(a)$ is of the form \cref{M21} the generator $G$ is of the form 
		\[
			G=
			\begin{pmatrix}
				-p & 0
				\\
				-\psi & p
			\end{pmatrix}
			,
		\]
		with some $\psi\in\bb R$. Now we obtain 
		\begin{align*}
			& \Psi(a)=\Psi\big(e^{-(-\log a)}\big)
			=\exp\Big(\smmatrix{-p}{0}{-\psi}{p}(-\log a)\Big)
			\\[3mm]
			&\mkern100mu =
			\begin{cases}
				\begin{pmatrix} 1 & 0 \\ \frac{\psi}{2p} & 1 \end{pmatrix}
				\begin{pmatrix} a^p & 0 \\ 0 & a^{-p} \end{pmatrix}
				\begin{pmatrix} 1 & 0 \\ -\frac{\psi}{2p} & 1 \end{pmatrix}
				&\text{if}\ p\neq 0
				,
				\\[6mm]
				\begin{pmatrix} 1 & 0 \\ \psi\log a & 1 \end{pmatrix}
				&\text{if}\ p=0
				,
			\end{cases}
		\end{align*}
		i.e., $\Psi(a)=\ms D_\psi(a)$. 
	\item We exploit \cref{M19} to show that $W(a,b,z)$ satisfies an integral equation and determine the Hamiltonian.
		Using that relation with $c=\frac 1b$ and evaluating derivatives w.r.t.\ $z$ at $z=0$ leads to 
		\[
			\forall 0<a\leq b\leq 1 \DP 
			\ms D_\psi(a)^{-1}\Big(\frac 1b W'(a,b,0)\Big)\ms D_\psi(b)
			=\ms D_\psi\big(\frac ab\big)^{-1}W'\big(\frac ab,1,0\big)
			.
		\]
		Since $\ms D_\psi(b)J\ms D_\psi(b)^T=J$, we obtain 
		\begin{equation}
		\label{M24}
			\forall 0<a<b\leq 1 \DP 
			\frac{W'(a,b,0)J}{b-a}=\ms D_\psi(b)\frac{W'(\frac ab,1,0)J}{1-\frac ab}\ms D_\psi(b)^T
			.
		\end{equation}
		The function $M:a\mapsto -W'(a,1,0)J$ takes nonpositive matrices as values, is nondecreasing in the matrix 
		sense, and $M(1)=0$. Thus its diagonal entries are nondecreasing and its off-diagonal entry is of bounded 
		variation. Since $W'(a,1,0)=W_E'(\alpha^{-1}(a),0,0)$, the function $M$ is also continuous. Hence, $M$ is 
		differentiable almost everywhere. In particular, there exists $b\in(0,1]$ such that the limit 
		\[
			\lim_{a\uparrow b}\frac{W'(a,b,0)J}{b-a}=\lim_{a\uparrow b}\frac{M(b)-M(a)}{b-a}
		\]
		exists (recall here the multiplicativity property of fundamental solutions). 
		Reading \cref{M24} from left to right, and remembering that $\ms D_\psi$ is continuous, 
		yields that the limit 
		\[
			P\DE\lim_{c\uparrow 1}\frac{W'(c,1,0)J}{1-c}
		\]
		exists. Clearly, $P$ is positive semidefinite. Reading \cref{M24} from right to left yields that $M$ is 
		everywhere differentiable and has the continuous derivative 
		\[
			\frac{dM}{da}(a)=\ms D_\psi(a)P\ms D_\psi(a)^T=H_{P,\psi}(a)\quad\text{for }a\in(0,1]
			.
		\]
		In particular, $M$ is absolutely continuous and can be written as the integral of its derivative. 
		Since $M$ is not constant, we must have $P\neq 0$. 

		Let us make a change of variable with the absolutely continuous 
		function $\alpha:(t_-,0]\to(0,1]$ from \cref{M18}. This gives, for each $a\in(0,1]$, 
		\begin{multline*}
			\int_{\alpha^{-1}(a)}^0 H_E(t) \DD t=W_E'(\alpha^{-1}(a),0,0)J=W'(a,1,0)J
			\\
			=M(1)-M(a)=\int_a^1 H_{P,\psi}(c) \DD c=\int_{\alpha^{-1}(a)}^0 H_{P,\psi}(\alpha(t))\alpha'(t) \DD t
			.
		\end{multline*}
		It follows that
		\begin{equation}
		\label{M26}
			H_E(t)=H_{P,\psi}(\alpha(t))\alpha'(t)\quad\text{for }t\in(t_-,0)\text{ a.e.}
			.
		\end{equation}
		Here is the point where we see that $\alpha^{-1}$ is absolutely continuous, since the above relation 
		shows $\alpha'(t)>0$ a.e. Thus, we may say that $H_E|_{(t_-,0)}$ is a 
		reparameterisation of $H_{P,\psi}$. Since $H_E|_{(t_-,0)}$ does not contain any indivisible intervals, also 
		$H_{P,\psi}|_{(0,1)}$ does not contain any such intervals. \Cref{M41} implies that $(P,\psi)\in\bb P_p$. 
	\item Using the same change of variable that led to \cref{M26} gives the integral equation for $W(a,b,z)$: 
		for all $0<a\leq b\leq 1$, we have
		\begin{align*}
			z\int_a^b W(a,c,z) &\, H_{P,\psi}(c) \DD c=
			z\int_{\alpha^{-1}(a)}^{\alpha^{-1}(b)} W(a,\alpha(t),z)H_{P,\psi}(\alpha(t))\alpha'(t)\DD t
			\\
			= &\, z\int_{\alpha^{-1}(a)}^{\alpha^{-1}(b)} W_E(\alpha^{-1}(a),t,z)H_E(t) \DD t
			\\
			= &\, W_E(\alpha^{-1}(a),\alpha^{-1}(b),z)J-J=W(a,b,z)J-J
			.
		\end{align*}
		Using this relation for $b=1$, multiplying from the left with $(a\odot_p A,a\odot_p B)\ms D_\psi(a)^{-1}$, 
		and remembering \cref{M33}, yields 
		\begin{multline*}
			\big(A(z),B(z)\big)J-\big([a\odot_p A](z),[a\odot_p B](z)\big)\ms D_\psi(a)^{-1}J
			\\
			=z\int_a^1 \big([c\odot_p A](z),[c\odot_p B](z)\big)\ms D_\psi(c)^{-1}H_{P,\psi}(c)\DD c
			.
		\end{multline*}
		\Cref{M53} now implies that \cref{M30} holds, i.e., $E=\widehat\Xi_p(P,\psi)$.
	\end{Elist}
\end{proof}

\begin{proof}[Proof of \Cref{M16}\,$(iii)$]
	For any two functions $E,\tilde E\in\HB$ with $E(0)=\tilde E(0)=1$ the spaces $\mc H(E)$ and $\mc H(\tilde E)$ are equal
	isometrically if and only if there exists $\gamma\in\bb R$ such that the structure Hamiltonian $H_{\tilde E}$ of 
	$\tilde E$ is a reparameterisation of $\smmatrix 10\gamma 1 H_E\smmatrix 1\gamma 01$ where $H_E$ is the structure
	Hamiltonian of $E$. Using the already proven statement \Cref{M16}\,$(i)$, this leads to 
	\begin{multline*}
		\mc H(\widehat\Xi_p(P,\psi))=\mc H(\widehat\Xi_p(\tilde P,\tilde\psi))\text{ isometrically}
		\quad\Longleftrightarrow\quad
		\\
		\exists\gamma\in\bb R\DP H_{\tilde P,\tilde\psi}|_{(0,1)}\text{ is a reparameterisation of }
		\smmatrix 10\gamma 1 H_{P,\psi}|_{(0,1)}\smmatrix 1\gamma 01
		.
	\end{multline*}
	The condition on the right means that there exists an absolutely continuous bijection 
	$\varphi\DF(0,1)\to(0,1)$ with $\varphi'>0$ a.e., such that 
	\begin{equation}
	\label{M61}
		\ms D_{\tilde\psi}(a)\tilde P\ms D_{\tilde\psi}(a)^T=
		\smmatrix 10\gamma 1\Big[\ms D_\psi P\ms D_\psi^T\Big]\big(\varphi(a)\big)\smmatrix 1\gamma 01\cdot\varphi'(a)
		.
	\end{equation}
	We write $P=\smmatrix{\kappa_1}{\kappa_3}{\kappa_3}{\kappa_2}$ and 
	$\tilde P=\smmatrix{\tilde\kappa_1}{\tilde\kappa_3}{\tilde\kappa_3}{\tilde\kappa_2}$. 
	Comparing the left upper entries, shows that \cref{M61} implies 
	\[
		\forall a\in(0,1)\DP a^{2p}\tilde\kappa_1=a^{2p}\kappa_1\cdot\varphi'(a)
		,
	\]
	and hence $\varphi=\Id_{(0,1)}$. We see that 
	\begin{multline}
	\label{M43}
		\mc H(\widehat\Xi_p(P,\psi))=\mc H(\widehat\Xi_p(\tilde P,\tilde\psi))\text{ isometrically}
		\quad\Longleftrightarrow\quad
		\\
		\exists\gamma\in\bb R\DQ\forall a\in(0,1)\DP 
		\ms D_{\tilde\psi}(a)\tilde P\ms D_{\tilde\psi}(a)^T=
		\smmatrix 10\gamma 1\ms D_\psi(a)P\ms D_\psi(a)^T\smmatrix 1\gamma 01
		.
	\end{multline}
	\begin{Ilist}
	\item \emph{Case $p=0$:}
		We plug the definition of $\ms D_\psi(a)$ and $\ms D_{\tilde\psi}(a)$ in \cref{M43}. 
		Remembering that $\kappa_1\neq 0$, this leads to 
		\begin{align*}
			& \mc H(\widehat\Xi_p(P,\psi))=\mc H(\widehat\Xi_p(\tilde P,\tilde\psi))\text{ isometrically}
			\\
			&\mkern50mu \Longleftrightarrow\quad
			\exists\gamma\in\bb R\DQ\forall a\in(0,1)\DP 
			\\
			&\mkern190mu
			\tilde P=
			\smmatrix 10{(\psi-\tilde\psi)\log a+\gamma}1\cdot P\cdot \smmatrix 1{(\psi-\tilde\psi)\log a+\gamma}01
			\\
			&\mkern50mu \Longleftrightarrow\quad 
			\exists\gamma\in\bb R\DQ\forall a\in(0,1)\DP 
			\tilde\kappa_1=\kappa_1\wedge\det\tilde P=\det P
			\\
			&\mkern265mu
			\wedge\tilde\kappa_3=\kappa_3+\kappa_1\big[(\psi-\tilde\psi)\log a+\gamma\big]
			\\
			&\mkern50mu \Longleftrightarrow\quad 
			\tilde\kappa_1=\kappa_1\wedge\det\tilde P=\det P\wedge\psi-\tilde\psi=0
			\\
			&\mkern50mu \Longleftrightarrow\quad 
			\tilde\kappa_1=\kappa_1\wedge\det\tilde P=\det P\wedge
			\psi-\tilde\psi=\frac{2p}{\kappa_1}(\kappa_3-\tilde\kappa_3)
			.
		\end{align*}
	\item \emph{Case $p\neq 0$:}
		In this case the computation is a bit more complicated; it is based on the formula
		\begin{equation}
		\label{M44}
			\begin{pmatrix} 1 & 0\\ \beta & 1\end{pmatrix}\cdot
			\begin{pmatrix} \alpha & 0 \\ 0 & \frac1\alpha \end{pmatrix}=
			\begin{pmatrix} \alpha & 0 \\ 0 & \frac1\alpha \end{pmatrix}\cdot
			\begin{pmatrix} 1 & 0\\ \alpha^2\beta & 1\end{pmatrix}
		\end{equation}
		which holds for all $\alpha,\beta\in\bb R$. 
		
		Plugging the formulae for $\ms D_\psi(a)$ and $\ms D_{\tilde\psi}(a)$, and using \cref{M44}, yields 
		\begin{align*}
			& \mc H(\widehat\Xi_p(P,\psi))=\mc H(\widehat\Xi_p(\tilde P,\tilde\psi))\text{ isometrically}
			\\
			&\mkern50mu \Longleftrightarrow\quad
			\exists\gamma\in\bb R\DQ\forall a\in(0,1)\DP 
			\\
			&\mkern75mu
			\smmatrix{a^p}00{a^{-p}}
			\smmatrix 10{-\frac{\tilde\psi}{2p}}1
			\cdot\tilde P\cdot 
			\smmatrix 1{-\frac{\tilde\psi}{2p}}01
			\smmatrix{a^p}00{a^{-p}}
			=
			\\
			&\mkern75mu
			\smmatrix{a^p}00{a^{-p}}
			\smmatrix 10{a^{2p}\big(\gamma+\frac{\psi-\tilde\psi}{2p}\big)-\frac{\psi}{2p}}1
			\cdot P\cdot 
			\smmatrix 1{a^{2p}\big(\gamma+\frac{\psi-\tilde\psi}{2p}\big)-\frac{\psi}{2p}}01
			\smmatrix{a^p}00{a^{-p}}
			\\
			&\mkern50mu \Longleftrightarrow\quad 
			\exists\gamma\in\bb R\DQ\forall a\in(0,1)\DP 
			\\
			&\mkern155mu
			\tilde P=
			\smmatrix 10{a^{2p}\big(\gamma+\frac{\psi-\tilde\psi}{2p}\big)-\frac{\psi-\tilde\psi}{2p}}1
			\cdot P\cdot 
			\smmatrix 1{a^{2p}\big(\gamma+\frac{\psi-\tilde\psi}{2p}\big)-\frac{\psi-\tilde\psi}{2p}}01
			\\
			&\mkern50mu \Longleftrightarrow\quad 
			\exists\gamma\in\bb R\DQ\forall a\in(0,1]\DP 
			\tilde\kappa_1=\kappa_1\wedge\det\tilde P=\det P
			\\
			&\mkern230mu
			\wedge\tilde\kappa_3=\kappa_3+
			\kappa_1\Big[a^{2p}\Big(\gamma+\frac{\psi-\tilde\psi}{2p}\Big)-\frac{\psi-\tilde\psi}{2p}\Big]
			\\
			&\mkern50mu \Longleftrightarrow\quad 
			\tilde\kappa_1=\kappa_1\wedge\det\tilde P=\det P\wedge
			\psi-\tilde\psi=\frac{2p}{\kappa_1}(\kappa_3-\tilde\kappa_3)
			.
		\end{align*}
	\end{Ilist}
\end{proof}

\noindent
Items $(i)$ and $(ii)$ of \Cref{M16} say that the map $(P,\psi)\mapsto\mc H(\widehat\Xi_p(P,\psi))$ is a surjection of the 
parameter space $\bb P_p$ onto the set of all homogeneous de~Branges spaces of order $\nu\DE p-\frac 12$, and in item $(iii)$ of
the theorem the equivalence relation corresponding to this map is described. We write $\approx$ for this relation, i.e.,
\begin{align*}
	(P,\psi)\approx(\tilde P,\tilde\psi)\DI &\, 
	\binom 10^*P\binom 10=\binom 10^*\tilde P\binom 10\ \wedge\ \det P=\det\tilde P
	\\
	&\ \wedge\ \psi-\tilde\psi=\frac{2p}{\binom 10^*P\binom 10}\Big[\binom 01^*P\binom 10-\binom 01^*\tilde P\binom 10\Big]
	.
\end{align*}
Naturally, we are interested in having at hand complete and minimal systems of representatives of our parameter space $\bb P_p$
modulo $\approx$. 

\begin{lemma}
\label{M42}
	Let $p>-\frac 12$. 
	\begin{Enumerate}
	\item The set 
		\begin{equation}
		\label{M45}
			\Big\{(P,\psi)\in\bb P_p \DS \binom 10^*P\binom 01=0\Big\}
		\end{equation}
		is a complete and minimal system of representatives of $\bb P_p$ modulo $\approx$. 
	\item If $p\neq 0$, then also 
		\begin{equation}
		\label{M46}
			\big\{(P,\psi)\in\bb P_p \DS \psi=0\big\}
		\end{equation}
		is a complete and minimal system of representatives of $\bb P_p$ modulo $\approx$. 
	\end{Enumerate}
\end{lemma}
\begin{proof}
	Let $(P,\psi)\in\bb P_p$ be given and write, as usual, $P=\smmatrix{\kappa_1}{\kappa_3}{\kappa_3}{\kappa_2}$. 
	We set
	\[
		\tilde P\DE \begin{pmatrix} \kappa_1 & 0 \\ 0 & \kappa_2-\frac{\kappa_3^2}{\kappa_1} \end{pmatrix}
		,\quad \tilde\psi\DE \psi-\frac{2p}{\kappa_1}\kappa_3
		.
	\]
	It is clear that $\tilde P\geq 0$, that $\det\tilde P=\det P$, $\tilde\kappa_1=\kappa_1$, and that 
	$\psi-\tilde\psi=\frac{2p}{\kappa_1}(\kappa_3-\tilde\kappa_3)$. 
	We certainly have $\binom 10\notin\ker\tilde P$. If $p\neq 0$ and $\tilde P\binom{-\tilde\psi}{2p}=0$, then 
	also $P\binom{-\psi}{2p}=0$ which is not possible. If $p=0$ and $\tilde\psi=0$, then also $\psi=0$ and hence 
	$\det P\neq 0$, which implies that $\det\tilde P\neq 0$. Hence, in all cases, $(\tilde P,\tilde\psi)\in\bb P_p$, and
	we see that $(\tilde P,\tilde\psi)\approx(P,\psi)$. If we have $(P,\psi),(\tilde P,\tilde\psi)\in\bb P_p$ 
	with $\kappa_3=\tilde\kappa_3=0$ which are in relation $\approx$, then the definition of $\approx$ immediately
	implies that $(P,\psi)=(\tilde P,\tilde\psi)$. We see that \cref{M45} is a complete and minimal system of representatives.

	Assume now that $p\neq 0$. Given $(P,\psi)\in\bb P_p$, we set 
	\[
		\tilde\kappa_1\DE\kappa_1,\ \tilde\kappa_3\DE\kappa_3-\frac{\kappa_1}{2p}\psi,\ 
		\tilde\kappa_2\DE\frac 1{\kappa_1}\big(\det P+\tilde\kappa_3^2\big),\quad 
		\tilde\psi\DE 0
		.
	\]
	Then, clearly, $\det\tilde P=\det P$, $\tilde P\geq 0$, and $\binom 10\notin\ker\tilde P$. 
	Assume towards a contradiction that $\tilde P\binom{-\tilde\psi}{2p}=0$. Then $\tilde\kappa_2=0$ and 
	$\det P=0$. The first implies that $\kappa_3=\frac{\kappa_1}{2p}\psi$ and the second that 
	$\kappa_2=\frac{\kappa_3^2}{\kappa_1}$. From this we see that $P\binom{-\psi}{2p}=0$, which is a 
	contradiction. Thus $(\tilde P,\tilde\psi)\in\bb P_p$, and of course the definition of $(\tilde P,\tilde\psi)$ 
	is made such that $(P,\psi)\approx(\tilde P,\tilde\psi)$. 
	If we have $(P,\psi),(\tilde P,\tilde\psi)\in\bb P_p$ with $\psi=\tilde\psi=0$ which are in relation $\approx$, 
	then immediately $(P,\psi)=(\tilde P,\tilde\psi)$. Thus, \cref{M46} is a complete and minimal system of representatives.
\end{proof}

\begin{Remark}
\label{M17}
	Let us discuss the case that $\nu=-\frac 12$ (i.e., $p=0$) in some more detail. 
	The de~Branges spaces which are homogeneous of order $-\frac 12$ are in one-to-one correspondence to
	the parameters 
	\[
		(\kappa_1,\kappa_2,\psi)\in
		\Big[(0,\infty)\times(0,\infty)\times\{0\}\Big]\cup
		\Big[(0,\infty)\times[0,\infty)\times\big(\bb R\setminus\{0\}\big)\Big]
		.
	\]
	Set $P\DE\smmatrix{\kappa_1}00{\kappa_2}$. If $\psi=0$, and hence $\kappa_2>0$, we find by inspecting the 
	formulae \cref{M123} and \cref{M124} that 
	\[
		\widehat\Xi_0(P,0)=
		\cos(\sqrt{\kappa_1\kappa_2}\cdot z)-i\sqrt{\frac{\kappa_1}{\kappa_2}}\sin(\sqrt{\kappa_1\kappa_2}\cdot z)
		.
	\]
	This family includes the Paley-Wiener spaces, namely for $\kappa_1=\kappa_2$. 
	If $\psi\neq 0$, the formulae are much more complicated and involve Kummer functions whose first argument is purely
	imaginary and nonzero, cf.\ \Cref{M36}. 

	In \cite{debranges:1962b,debranges:1968} the false argument occurs that \emph{for every} $\nu>-\frac 12$ the totality
	of homogeneous spaces is obtained with parameter $\psi=0$. This is true for $p\neq 0$, but for $p=0$ 
	the whole family of spaces occurring from 
	parameters $(\kappa_1,\kappa_2,\psi)\in(0,\infty)\times[0,\infty)\times(\bb R\setminus\{0\})$ was lost.
\end{Remark}

%\newpage
%
%
%
\section[\textcolor{ForestGreen}{Measures induced by homogeneous spaces}]{Measures induced by homogeneous spaces}
%ForestGreen / Dandelion / BrickRed / Fuchsia
%
%

Let $\nu>-1$ and set, as usual, $p\DE\nu+\frac 12$. If $\mc H(E)$ is a de~Branges space which is homogeneous of order $\nu$, 
\Cref{M7}\,$(i)$ shows that 
\begin{equation}
\label{M38}
	\big\{\mc H(a\odot_p E) \DS a\in(0,\infty)\big\}
\end{equation}
is a chain of de~Branges spaces with isometric inclusions $\mc H(a\odot_p E)\subseteq\mc H(b\odot_p E)$ when $a\leq b$. 

The chain \cref{M38} has no gaps or jumps, i.e., 
\[
	\forall b\in(0,\infty)\DP 
	\bigcap_{a>b}\mc H(a\odot_p E)=\mc H(b\odot_p E)\wedge
	\Clos\bigcup_{a<b}\mc H(a\odot_p E)=\mc H(b\odot_p E)
	,
\]
cf.\ \Cref{M14}. Moreover, due to the kernel relation \cref{M4}, it is exhaustive in the sense that 
\[
	\lim_{a\downarrow 0}K_{a\odot_p E}(z,w)=0,\quad \lim_{a\uparrow\infty}K_{a\odot_p E}(0,0)=\infty
	.
\]
These properties show that $\{\mc H(a\odot_p E)\DS a>0\}$ is an unbounded chain in the sense of \Cref{M22}. 
By \Cref{M23} it determines a positive Borel measure on the real line.

\begin{Definition}
\label{M40}
	Let $\nu>-1$ and let $\mc H$ be a de~Branges space which is homogeneous of order $\nu$, and choose $E\in\HB$ such that
	$\mc H=\mc H(E)$. 
	Then we denote by $\mu_{\mc H}$ the unique positive Borel measure on the real line, such
	that $L^2(\mu_{\mc H})$ contains $\bigcup_{a\in(0,\infty)}\mc H(a\odot_p E)$ isometrically as a dense subspace. 

	Morever, let us introduce the following abbreviation: 
	given $(P,\psi)\in\bb P_p$, we denote $\mu_{P,\psi}\DE\mu_{\mc H(\hat\Xi_p(P,\psi))}$.
\end{Definition}

\noindent
The measures occuring as $\mu_{\mc H}$ for some homogeneous space have a very particular form. 
We formulate this in terms of the parameter class $\bb P_p$. 

\begin{Theorem}
\label{M37}
	Let $\nu>-1$ and set $p\DE\nu+\frac 12$. Moreover, let $\lambda$ be the Lebesgue measure on $\bb R$. 
	Then the following statements hold. 
	\begin{Enumerate}
	\item Let $(P,\psi)\in\bb P_p$ and write $P=\smmatrix{\kappa_1}{\kappa_3}{\kappa_3}{\kappa_2}$ and 
		$\sigma\DE2p\kappa_3-\psi\kappa_1$. Moreover, let $\kappa$ be the nonnegative square root of $\det P$. 
		Then $\mu_{P,\psi}\ll\lambda$ and 
		\begin{equation}
		\label{M119}
			\frac{d\mu_{P,\psi}}{\DD\lambda}(x)=
			\begin{cases}
				\mu_+(P,\psi)\cdot |x|^{2p} \CAS x>0
				,
				\\[1mm]
				\mu_-(P,\psi)\cdot |x|^{2p} \CAS x<0
				,
			\end{cases}
		\end{equation}
		where 
		\begin{align}
		\label{M130}
			\mu_+(P,\psi)\DE &\, 
			\begin{cases}
				\frac 1\pi
				\frac{2^{2p}\kappa^{2p+1}|\Gamma(\frac{\sigma}{2i\kappa}+p+1)|^2}{\kappa_1\Gamma(2p+1)^2}
				\cdot e^{\pi\frac{\sigma}{2\kappa}}
				\CAS \det P\neq 0
				,
				\\
				\frac{\sigma^{2p+1}}{\kappa_1\Gamma(2p+1)^2}
				\CAS \det P=0,\sigma>0
				,
				\\
				0
				\CAS \det P=0,\sigma<0
				,
			\end{cases}
			\\[1mm]
		\label{M131}
			\mu_-(P,\psi)\DE &\, 
			\begin{cases}
				\frac 1\pi
				\frac{2^{2p}\kappa^{2p+1}|\Gamma(\frac{\sigma}{2i\kappa}+p+1)|^2}{\kappa_1\Gamma(2p+1)^2}
				\cdot e^{-\pi\frac{\sigma}{2\kappa}}
				\CAS \det P\neq 0
				,
				\\
				0
				\CAS \det P=0,\sigma>0
				,
				\\
				\frac{|\sigma|^{2p+1}}{\kappa_1\Gamma(2p+1)^2}
				\CAS \det P=0,\sigma<0
				.
			\end{cases}
		\end{align}
	\item Let $(\mu_+,\mu_-)\in[0,\infty)^2\setminus\{(0,0)\}$. Then there exists $(P,\psi)\in\bb P_p$ such that 
		$\mu_+=\mu_+(P,\psi)$ and $\mu_-=\mu_-(P,\psi)$. 
	\item Let $(P,\psi),(\tilde P,\tilde\psi)\in\bb P_p$, and write $P=\smmatrix{\kappa_1}{\kappa_3}{\kappa_3}{\kappa_2}$
		and $P=\smmatrix{\tilde\kappa_1}{\tilde\kappa_3}{\tilde\kappa_3}{\tilde\kappa_2}$.
		Then $\mu_{P,\psi}=\mu_{\tilde P,\tilde\psi}$, if and only if
		\begin{Ilist}
		\item $\kappa_1^{-\frac 2{1+2p}}\det P=\tilde\kappa_1^{-\frac 2{1+2p}}\det\tilde P$, 
		\item $\kappa_1^{\frac{2p}{1+2p}}\psi-\tilde\kappa_1^{\frac{2p}{1+2p}}\tilde\psi=
			2p\Big(\kappa_1^{-\frac 1{1+2p}}\kappa_3-\tilde\kappa_1^{-\frac 1{1+2p}}\tilde\kappa_3\Big)$
			.
		\end{Ilist}
	\end{Enumerate}
\end{Theorem}

\noindent
Concerning item $(i)$ of the theorem, in \cite{debranges:1962b,debranges:1968} the following less precise statement is shown.

\begin{proposition}
\label{M39}
	Let $\nu>-1$, set $p\DE\nu+\frac 12$, and let $(P,\psi)\in\bb P_p$. Denote by $\lambda$ the Lebesgue measure on $\bb R$. 
	Then $\mu_{P,\psi}\ll\lambda$ and its derivative is of the form \cref{M119} with some numbers 
	$\mu_\pm(P,\psi)\geq 0$. 
\end{proposition}
\begin{proof}
	Let $c>0$. Then the map $F(x)\mapsto[c\odot_{p+\frac 12} F](x)$ is an isometric bijection of 
	$\bigcup_{a>0}\mc H(a\odot_p E)$ onto itself. Thus it has an extension to a unitary operator of $L^2(\mu_{P,\psi})$ 
	onto itself. Since $L^2$-convergence implies pointwise a.e.\ convergence of a subsequence, this extension acts again as 
	$f(x)\mapsto[c\odot_{p+\frac 12} f](x)$ (for $x\in\bb R$ a.e.). 

	We have, for every $c>0$, 
	\[
		\mu_{P,\psi}\big((0,c)\big)=\|\mathds{1}_{(0,c)}\|_{L^2(\mu_{P,\psi})}^2
		=\|\underbrace{c\odot_{p+\frac 12}\mathds{1}_{(0,c)}}_{=c^{p+\frac 12}\mathds{1}_{(0,1)}}\|_{L^2(\mu_{P,\psi})}^2
		=c^{2p+1}\mu_{P,\psi}\big((0,1)\big)
		.
	\]
	Analogously, we find 
	\[
		\mu_{P,\psi}\big((-c,0)\big)=c^{2p+1}\mu_{P,\psi}\big((-1,0)\big)
		\text{ and }
		\mu_{P,\psi}\big((-c,c)\big)=c^{2p+1}\mu_{P,\psi}\big((-1,1)\big)
		.
	\]
	The first relation shows that $\mathds{1}_{(0,\infty)}\DD\mu_{P,\psi}$ is absolutely continuous w.r.t.\ $\lambda$, and 
	\[
		\frac{\mathds{1}_{(0,\infty)}(x)\DD\mu_{P,\psi}(x)}{\DD\lambda(x)}=
		\mathds{1}_{(0,\infty)}(x)(2p+1)x^{2p}\mu_{P,\psi}\big((0,1)\big)
		,
	\]
	the second that $\mathds{1}_{(-\infty,0)}\DD\mu_{P,\psi}$ is absolutely continuous w.r.t.\ $\lambda$, and 
	\[
		\frac{\mathds{1}_{(-\infty,0)}(x)\DD\mu_{P,\psi}(x)}{\DD\lambda(x)}=
		\mathds{1}_{(-\infty,0)}(x)(2p+1)|x|^{2p}\mu_{P,\psi}\big((-1,0)\big)
		,
	\]
	and letting $c\downarrow 0$ in the third relation yields 
	\[
		\mu_{P,\psi}(\{0\})=\lim_{c\downarrow 0}\mu_{P,\psi}\big((-c,c)\big)=0
		.
	\]
\end{proof}

\noindent
The explicit formulae \cref{M130} and \cref{M131} for $\mu_\pm(P,\psi)$ are stated without a proof in 
\cite[p.210]{debranges:1962b}. Obtaining those values requires knowing the explicit formulae for $\hat\Xi_p(P,\psi)$ in terms 
of Kummer functions. The argument rests on the following fact about asymptotics. 

\begin{lemma}
\label{M132}
	Let $\delta\in\bb R$ and $p\in\bb R\setminus(-\frac 12\bb N)$. Then there exist bounded functions 
	$R_\pm\DF [1,\infty)\to\bb R$, such that 
	\begin{align*}
		& \Big|\frac 12 M(i\delta+p,2p+1,-iy)+\frac 12 M(i\delta+p+1,2p+1,-iy)
		\\
		&\mkern210mu
		-\frac i{2(2p+1)}yM(i\delta+p+1,2p+2,-iy)\Big|
		\\
		& =\frac{\Gamma(2p+1)}{|\Gamma(i\delta+p+1)|} e^{\pm\frac\pi 2\delta}\cdot\frac 1{|y|^p}
		\cdot\Big(1+\frac{R_\pm(|y|)}{|y|}\Big)
		,
	\end{align*}
	where the ``$+$''-sign holds if $y\geq 1$ and the ``$-$''-sign when $y\leq -1$. 
\end{lemma}
\begin{proof}
	We use the asymptotic expansion of the Kummer function given in \cite[\phantom{}13.5.1]{abramowitz.stegun:1964}, see also 
	\cite[\phantom{}13.7.2]{nist:2010}. For a purely imaginary argument $z=iy$, this reads as (for parameters $a\in\bb C$ and 
	$b\in\bb C\setminus(-\bb N_0)$)
	\begin{align*}
		M(a,b,iy)= &\, 
		\frac{\Gamma(b)}{\Gamma(b-a)}e^{\pm i\frac\pi 2a}e^{-i(\Im a)\log|y|}\cdot |y|^{-\Re a}
		\cdot\big(1+\BigO(\smfrac 1{|y|})\big)
		\\
		+ &\,
		\frac{\Gamma(b)}{\Gamma(a)}e^{\pm i\frac\pi 2(a-b)}e^{i(y+(\Im(a-b))\log|y|)}\cdot |y|^{\Re(a-b)}
		\cdot\big(1+\BigO(\smfrac 1{|y|})\big)
		,
	\end{align*}
	where the ``$+$''-sign holds if $y>0$ and the ``$-$''-sign if $y<0$, and where the $\BigO(\frac 1{|y|})$ 
	is understood for $|y|\to\infty$. 

	Let $\gamma\in\bb C$, then these formulae give 
	\begin{align*}
		& \frac 12 M(i\delta+p,2p+1,-iy)+\frac 12 M(i\delta+p+1,2p+1,-iy)
		\\
		&\mkern205mu
		-\gamma yM(i\delta+p+1,2p+2,-iy)
		\\
		& =\pm\smfrac{\Gamma(2p+2)}{\Gamma(-i\delta+p+1)} e^{\mp i\frac\pi 2(i\delta+p+1)}e^{-i\delta\log|y|}
		\cdot\frac 1{|y|^p}\cdot\big[\smfrac i{2(2p+1)}-\gamma\big]
		\cdot\big(1+\BigO(\smfrac 1{|y|})\big)
		\\
		& \pm\smfrac{\Gamma(2p+2)}{\Gamma(i\delta+p+1)} e^{\pm i\frac\pi 2(-i\delta+p+1)}e^{i(y+\delta\log|y|)}
		\cdot\frac 1{|y|^p}\cdot\big[\smfrac{-i}{2(2p+1)}-\gamma\big]
		\cdot\big(1+\BigO(\smfrac 1{|y|})\big)
	\end{align*}
	where the upper sign holds if $y>0$ and the lower sign if $y<0$.

	Using $\gamma\DE\frac i{2(2p+1)}$ in this formula yields the assertion of the lemma. 
\end{proof}

\begin{proof}[Proof of \Cref{M37}\,$(i)$; Case $\det P\neq 0$]
	As usual we write $\Xi_p(P,\psi)=(A,B)$ and $\hat\Xi_p(P,\psi)=E$. 
	\begin{Elist}
	\item For $a>0$ and $\tau$ in the open upper half-plane $\bb C_+$ consider the function 
		\[
			q_{a,\tau}(z)\DE\frac{[a\odot_p A](z)\tau+[a\odot_p B](z)}{-[a\odot_p B](z)\tau+[a\odot_p A](z)}
			,\quad z\in\bb C_+
			.
		\]
		Then $q_{a,\tau}$ is a Nevanlinna function, as computing the Nevanlinna kernel shows. 
		Since $A$ and $B$ have no common real zeroes and $\Im\tau>0$, 
		it has an analytic continuation to some domain containing the closed upper half-plane $\bb C_+\cup\bb R$. 
		Hence, the measure $\nu_{a,\tau}$ in the Nevanlinna integral representation of $q_{a,\tau}$ is absolutely
		continuous w.r.t.\ the Lebesgue measure and 
		\[
			\frac{\DD\nu_{a,\tau}}{\DD\lambda}(x)=\frac 1\pi\Im q_{a,\tau}(x),\quad x\in\bb R\text{ a.e.}
		\]
		The imaginary part of $q_{a,\tau}$ computes as 
		\[
			\Im q_{a,\tau}(x)=\frac{(\Im\tau)\cdot|[a\odot_p E](x)|^2}{|-[a\odot_p B](x)\tau+[a\odot_p A](x)|^2}
			.
		\]
		By \cite[Theorem~32]{debranges:1968} the space $\mc H(a\odot_p E)$ is contained contractively in 
		$L^2(\mu_{a,\tau})$ where $\mu_{a,\tau}$ is the measure which is mutually absolutely continuous with 
		$\nu_{a,\tau}$ and has derivative
		\[
			\frac{\DD\mu_{a,\tau}}{\DD\nu_{a,\tau}}(x)=\frac{1}{|[a\odot_p E](x)|^2}
			.
		\]
		Since the chain \cref{M38} contains no one-dimensional gaps, this inclusion is actually always isometric. 

		Let us note explicitly that from the above $\DD\mu_{a,\tau}\ll\DD\lambda$ and 
		\[
			\frac{\DD\mu_{a,\tau}}{\DD\lambda}(x)=
			\frac 1\pi\frac{\Im\tau}{|-[a\odot_p B](x)\tau+[a\odot_p A](x)|^2}
			.
		\]
	\item We choose $\tau$ such that \Cref{M132} and \Cref{M36} is applicable: 
		\[
			\tau\DE -\frac{\kappa_3}{\kappa_1}+i\frac\kappa{\kappa_1}
			.
		\]
		Then \Cref{M132} yields 
		\[
			\big|-[a\odot_p B](x)\tau+[a\odot_p A](x)\big|=
			\frac{\Gamma(2p+1)}{|\Gamma(-\frac{\sigma}{2i\kappa}+p+1)|}
			\cdot e^{\mp\frac\pi 2\frac{\sigma}{2\kappa}}\cdot\frac 1{|2\kappa x|^p}
			\cdot\big(1+\smfrac{R_\pm(|ax|)}{|ax|}\big)
			,
		\]
		and in turn 
		\begin{align*}
			\frac{\DD\mu_{a,\tau}}{\DD\lambda}(x)= &\, 
			\frac 1\pi\frac{\frac{\kappa}{\kappa_1}}{
			\frac{\Gamma(2p+1)^2}{|\Gamma(-\frac{\sigma}{2i\kappa}+p+1)|^2}
			\cdot e^{\mp\pi\frac{\sigma}{2\kappa}}\cdot\frac 1{|2\kappa x|^{2p}}}
			\cdot\big(1+\smfrac{R_\pm(|ax|)}{|ax|}\big)
			\\[2mm]
			= &\, 
			\frac 1\pi\frac{2^{2p}\kappa^{2p+1}|\Gamma(\frac{\sigma}{2i\kappa}+p+1)|^2}{\kappa_1\Gamma(2p+1)^2}
			\cdot e^{\pm\pi\frac{\sigma}{2\kappa}}\cdot|x|^{2p}
			\cdot\big(1+\smfrac{R_\pm(|ax|)}{|ax|}\big)
			\\[2mm]
			= &\, 
			\begin{cases}
				\mu_+(P,\psi)\cdot x^{2p}\cdot\big(1+\smfrac{R_+(|ax|)}{|ax|}\big)
				\CAS x>0
				,
				\\[2mm]
				\mu_-(P,\psi)\cdot|x|^{2p}\cdot\big(1+\smfrac{R_-(|ax|)}{|ax|}\big)
				\CAS x<0
				.
			\end{cases}
		\end{align*}
	\item By \Cref{M23} we have $\lim_{a\to\infty}\mu_{a,\tau}=\mu_{P,\psi}$ vaguely. Since 
		\[
			\lim_{a\to\infty}\frac{\DD\mu_{a,\tau}}{\DD\lambda}(x)=
			\begin{cases}
				\mu_+(P,\psi)\cdot x^{2p}
				\CAS x>0
				,
				\\
				\mu_-(P,\psi)\cdot|x|^{2p}
				\CAS x<0
				,
			\end{cases}
		\]
		and the convergence is uniform on compact subsets of $\bb R\setminus\{0\}$, we obtain that 
		$\mathds{1}_{\bb R\setminus\{0\}}\cdot\mu_{P,\psi}\ll\lambda$ and 
		\[
			\frac{\DD(\mathds{1}_{\bb R\setminus\{0\}}\cdot\mu_{P,\psi})}{\DD\lambda}(x)=
			\begin{cases}
				\mu_+(P,\psi)\cdot x^{2p}
				\CAS x>0
				,
				\\
				\mu_-(P,\psi)\cdot|x|^{2p}
				\CAS x<0
				.
			\end{cases}
		\]
		As we saw in \Cref{M39}, $\mu_{P,\psi}$ has no point mass at $0$, and \cref{M119} follows. 
	\end{Elist}
\end{proof}

\begin{proof}[Proof of \Cref{M37}\,$(i)$; Case $\det P=0$]
	We use a continuity argument. Note that, clearly, the set of all $(P,\psi)\in\bb P_p$ with $\det P\neq 0$ is 
	dense in $\bb P_p$. The function $\hat\Xi_p(P,\psi)$, as well as the Hamiltonian $H_{P,\psi}$, depend continuously on 
	$(P,\psi)$ w.r.t.\ locally uniform convergence (on $\bb C$ for the first, and on $(0,\infty)$ for the latter).
	Hence, the measure $\mu_{P,\psi}$ depends continuously on $(P,\psi)$ w.r.t.\ vague convergence of measures,
	cf.\ \Cref{M62}.

	Let us show that $\det P=0$ implies $\sigma\neq 0$. Consider the case that $p\neq 0$. 
	We can write $\sigma=(\kappa_1,\kappa_3)\cdot\binom{-\psi}{2p}$, and since the rows of $P$ are linearly
	dependent, $\sigma=0$ implies that $\binom{-\psi}{2p}\in\ker P$. This is excluded by the definition of $\bb P_p$. 
	If $p=0$, then we must have $\psi\neq 0$ by the definition of $\bb P_p$, and also in this case it follows that 
	$\sigma\neq 0$. 

	In order to prove the second and third lines in \cref{M130} and \cref{M131}, we thus have to evaluate the limit of 
	the constants in the respective first lines when $\frac{\sigma}{\kappa}\to\pm\infty$. This is easy using the 
	relation
	\[
		|\Gamma(x+iy)|\sim\sqrt{2\pi}\cdot |y|^{x-\frac 12}e^{-\pi\frac{|y|}2}\quad\text{for }y\to\pm\infty
		, 
	\]
	cf.\ \cite[\phantom{}5.11.9]{nist:2010}. Namely, we obtain that (for $\kappa\to 0$) 
	\[
		\frac{2^{2p}\kappa^{2p+1}|\Gamma(\frac{\sigma}{2i\kappa}+p+1)|^2)}{\kappa_1\Gamma(2p+1)^2} 
		\cdot e^{\pi\frac{\sigma}{2\kappa}}
		\sim \frac{\pi|\sigma|^{2p+1}}{\kappa_1\Gamma(2p+1)^2} e^{\frac\pi{2\kappa}(\sigma-|\sigma|)}
		,
	\]
	and this yields the second and third lines in \cref{M130}. The relations in \cref{M131} follow analogously. 
\end{proof}

\begin{proof}[Proof of \Cref{M37}\,$(ii)$]
	Let $(\mu_+,\mu_-)\in[0,\infty)^2\setminus\{(0,0)\}$ be given. If $\mu_-=0$ or $\mu_+=0$, we use 
	$P\DE\smmatrix 1000$ and 
	\[
		\psi\DE-\Big(\mu_+\Gamma(2p+1)^2\Big)^{\frac 1{2p+1}}\quad\text{or}\quad
		\psi\DE\Big(\mu_-\Gamma(2p+1)^2\Big)^{\frac 1{2p+1}}
		,
	\]
	respectively. 

	If $\mu_+,\mu_->0$, we set 
	\[
		\kappa_1\DE 1,\ \kappa_2\DE\bigg(
		\frac{\pi\Gamma(2p+1)^2\sqrt{\mu_+\mu_-}}{2^{2p}
		\big|\Gamma\big(\frac i{2\pi}\log\frac{\mu_-}{\mu_+}+p+1\big)\big|^2}
		\bigg)^{\frac 1{p+\frac 12}}
		,\ \kappa_3\DE 0,\ \psi\DE\sqrt{\kappa_2}\cdot\frac 1\pi\log\frac{\mu_-}{\mu_+}
		.
	\]
	Then $(P,\psi)\in\bb P_p$, and plugging in the definitions yields that $\mu_\pm(P,\psi)=\mu_\pm$. 
\end{proof}

\begin{proof}[Proof of \Cref{M37}\,$(iii)$]
	Set $E\DE\hat\Xi_p(P,\psi)$ and $\tilde E\DE\hat\Xi_p(\tilde P,\tilde\psi)$. 
	\begin{Elist}
	\item We show that
		\[
			\mu_{P,\psi}=\mu_{\tilde P,\tilde\psi}
			\quad\Longleftrightarrow\quad
			\exists c\in\bb R\DP \mc H(\tilde E)=\mc H(c\odot_p E)\text{ isometrically}
			.
		\]
		The implication ``$\Leftarrow$'' follows since $\mc H(\tilde E)=\mc H(c\odot_p E)$ implies 
		\[
			\big\{\mc H(a\odot_p\tilde E)\DS a>0\big\}=\big\{\mc H(a\odot_p E)\DS a>0\big\}
			,
		\]
		and thus $L^2(\mu_{P,\psi})=L^2(\mu_{\tilde P,\tilde\psi})$. To show the implication ``$\Rightarrow$'', 
		assume that $\mu_{P,\psi}=\mu_{\tilde P,\tilde\psi}$. 
		The functions $E$ and $\tilde E$ are of bounded type in the upper half-plane and have no 
		real zeroes, cf.\ \Cref{M104} and \Cref{M51}. Hence, the Ordering Theorem \cite[Theorem~35]{debranges:1968} 
		applies and yields 
		\[
			\forall a>0\DP \mc H(\tilde E)\subseteq\mc H(a\odot_p E)\vee\mc H(a\odot_p E)\subseteq\mc H(\tilde E)
			.
		\]
		\Cref{M14} implies $\mc H(\tilde E)=\mc H(c\odot_p E)$ with 
		$c\DE \big(\frac{K_{\tilde E}(0,0)}{K_E(0,0)}\big)^{\frac 1{2p+1}}$.
	\item Given $c>0$, we construct $(P_c,\psi_c)\in\bb P_p$ such that $\mc H(\hat\Xi_p(P_c,\psi_c))=\mc H(c\odot_pE)$. 
		Set 
		\[
			(A_c,B_c)\DE(c\odot_pA,c\odot_pB)\begin{pmatrix} c^{-p} & 0\\ 0 & c^p\end{pmatrix}
			,
		\]
		then $E_c\DE A_c-iB_c$ satisfies $\mc H(E_c)=\mc H(c\odot_pE)$ by \Cref{M99}. 
		Write $A(z)=\sum_{n=0}^\infty \alpha_nz^n$ and $B(z)=\sum_{n=0}^\infty \beta_nz^n$, then 
		\[
			A(cz)=\sum_{n=0}^\infty \alpha_nc^n\cdot z^n,\quad c^{2p}B(cz)=\sum_{n=0}^\infty \beta_nc^{2p+n}\cdot z^n
			.
		\]
		A short computation using \cref{M44} shows that 
		\begin{multline*}
			(\alpha_{n+1}c^{n+1},\beta_{n+1}c^{2p+n+1})=-\frac 1{(n+1)(2p+n+1)}(\alpha_nc^n,\beta_nc^{2p+n})
			\\
			\cdot
			\smmatrix{c^{\frac 12+p}}00{c^{\frac 12-p}}P\smmatrix{c^{\frac 12+p}}00{c^{\frac 12-p}}
			\cdot J\smmatrix{2p+n+1}{0}{c^{-2p}\psi}{n+1}
		\end{multline*}
		for all $n\in\bb N_0$. 

		The pair 
		\begin{equation}
		\label{M48}
			(P_c,\psi_c)\DE
			\bigg(
			\smmatrix{c^{\frac 12+p}}00{c^{\frac 12-p}}P\smmatrix{c^{\frac 12+p}}00{c^{\frac 12-p}},c^{-2p}\psi
			\bigg)
		\end{equation}
		belongs to $\bb P_p$, and by the above computation $\hat\Xi_p(P_c,\psi_c)=E_c$. 
	\item We use \Cref{M16}\,$(iii)$ to complete the proof. This theorem tells us that 
		$\mc H(\hat\Xi_p(P_c,\psi_c))=\mc H(\hat\Xi_p(\tilde P,\tilde\psi))$ isometrically, if and only if 
		\begin{Ilist}
		\item $c^{1+2p}\kappa_1=\tilde\kappa_1$,
		\item $c^2\det P=\det\tilde P$,
		\item $c^{-2p}\psi-\tilde\psi=\frac{2p}{c^{1+2p}\kappa_1}\big(c\kappa_3-\tilde\kappa_3\big)$.
		\end{Ilist}
		The first relation gives 
		\[
			c=\Big(\frac{\tilde\kappa_1}{\kappa_1}\Big)^{\frac 1{1+2p}}
			.
		\]
		Plugging this into the second and third relations, leads to the stated formulae. 
	\end{Elist}
\end{proof}

\noindent
Items $(i)$ and $(ii)$ of \Cref{M37} say that the map $(P,\psi)\mapsto\mu_{P,\psi}$ is a surjection of the 
parameter space $\bb P_p$ onto the set of all measures $\mu\ll\lambda$ whose derivative is of the form 
\begin{equation}
\label{M111}
	\frac{\DD\mu}{\DD\lambda}(x)=
	\begin{cases}
		\mu_+\cdot x^{2p} \CAS x>0
		,
		\\
		\mu_-\cdot|x|^{2p} \CAS x<0
		,
	\end{cases}
\end{equation}
with some $(\mu_+,\mu_-)\in[0,\infty)^2\setminus\{(0,0)\}$. 
In item $(iii)$ of the theorem the equivalence relation corresponding to this map is described 
(we write $\simeq$ for this relation). 
In the next lemma, we provide a complete and minimal system of representatives of our parameter space $\bb P_p$ modulo $\simeq$. 

\begin{lemma}
\label{M64}
	Let $p>-\frac 12$. 
	\begin{Enumerate}
	\item The set 
		\begin{equation}
		\label{M84}
			\Big\{(P,\psi)\in\bb P_p \DS \binom 10^*P\binom 01=0\wedge\binom 10^*P\binom 10=1\Big\}
		\end{equation}
		is a complete and minimal system of representatives of $\bb P_p$ modulo $\simeq$. 
	\item If $p\neq 0$, then also 
		\begin{equation}
		\label{M85}
			\big\{(P,\psi)\in\bb P_p \DS \psi=0\wedge\binom 10^*P\binom 10=1\big\}
		\end{equation}
		is a complete and minimal system of representatives of $\bb P_p$ modulo $\simeq$. 
	\end{Enumerate}
\end{lemma}
\begin{proof}
	Let $(P,\psi)\in\bb P_p$ be given. Set $(\tilde P,\tilde\psi)\DE (P_c,\psi_c)$, cf.\ \cref{M48}, 
	with $c\DE\kappa_1^{-\frac 1{2p+1}}$. Then $(P,\psi)\simeq(\tilde P,\tilde\psi)$ and $\tilde\kappa_1=1$. 
	By the definition in \cref{M48} it is clear that $\kappa_3=0$ implies $\tilde\kappa_3=0$ and $\psi=0$ implies 
	$\tilde\psi=0$. 

	Since $\simeq\,\supseteq\,\approx$, it follows from \Cref{M42} that we can always reduce modulo $\simeq$ to an element of 
	the form written in \cref{M84} or \cref{M85}, respectively. 
	Assume we have two element $(P,\psi),(\tilde P,\tilde\psi)\in\bb P_p$ with $(P,\psi)\simeq(\tilde P,\tilde\psi)$ and 
	$\kappa_1=\tilde\kappa_1=1$ and $\kappa_3=\tilde\kappa_3=0$ (or $\psi=\tilde\psi=0$). Then 
	$\kappa_2=\tilde\kappa_2$ from the first relation in \Cref{M37}\,$(iii)$ and $\psi=\tilde\psi$ from the second 
	(or $\kappa_3=\tilde\kappa_3$ from the second and then $\kappa_2=\tilde\kappa_2$ from the first, respectively). 
\end{proof}

\noindent
Combining \Cref{M37} with \Cref{M36} we can directly connect measures \cref{M119} with their corresponding chains \cref{M38}. 

\begin{Corollary}
\label{M90}
	Let $p>-\frac 12$ and $(\mu_+,\mu_-)\in[0,\infty)^2\setminus\{(0,0)\}$. Let $\mu$ be the measure with $\mu\ll\lambda$ 
	and derivative \cref{M111}. We define functions $A,B$ by distinguishing cases.
	\begin{Enumerate}
	\item Assume that $\mu_+,\mu_->0$. Define 
		\begin{align*}
			& \alpha\DE \frac i{2\pi}\log\frac{\mu_-}{\mu_+}+p,\quad
			\kappa\DE\frac 12
			\bigg(\frac{\pi 2\Gamma(2p+2)^2\sqrt{\mu_+\mu_-}}{(2p+1)|\Gamma(\alpha+1)|^2}\bigg)^{\frac 1{2p+1}}
			,
			\\[1mm]
			& A(z)\DE e^{i\kappa z}\frac 12\Big[M(\alpha,2p+1,-2i\kappa z)+M(\alpha+1,2p+1,-2i\kappa z)\Big]
			,
			\\
			& B(z)\DE z\cdot e^{i\kappa z}M(\alpha+1,2p+2,-2i\kappa z)
			.
		\end{align*}
		Then $\mu_{\mc H(E)}=\mu$.
	\item Assume that $\mu_+=0$ or $\mu_-=0$. Define
		\begin{align*}
			& \sigma\DE
			\begin{cases}
				\Big(\frac{\mu_+}{2p+1}\Gamma(2p+2)^2\Big)^{\frac 1{2p+1}} \CAS \mu_+>0
				,
				\\
				-\Big(\frac{\mu_-}{2p+1}\Gamma(2p+2)^2\Big)^{\frac 1{2p+1}} \CAS \mu_->0
				,
			\end{cases}
			\\[2mm]
			& A(z)\DE \prescript{}{0}{F}_1(2p+1,-\sigma z),\quad B(z)= z\cdot\prescript{}{0}{F}_1(2p+2,-\sigma z)
			.
		\end{align*}
	\end{Enumerate}
	Then $\mu_{\mc H(E)}=\mu$.
\end{Corollary}

\noindent
Of course there are many choices for the function $E$ generating the chain of $\mu$. The choice in this corollary is made 
in such a way that $E(0)=1$ and $K_E(0,0)=1$. 

\begin{proof}[Proof of \Cref{M90}]
	In the proof of \Cref{M37}\,$(ii)$ we have already exhibited a pair $(P,\psi)$ such that the corresponding measure is
	equal to $\mu$. \Cref{M37}\,$(iii)$ allows us to modify this pair; and we use this freedom to obtain that 
	$K_E(0,0)=1$. 
	\begin{Ilist}
	\item If $\mu_+,\mu_->0$, set
		\begin{align*}
			& \kappa_1\DE 2p+1,\quad 
			\kappa_2\DE \Big(\frac 1{2p+1}\Big)^{1-\frac 1{p+\frac 12}} \bigg(
			\frac{\pi\Gamma(2p+1)^2\sqrt{\mu_+\mu_-}}{2^{2p}|\Gamma(\frac i{2\pi}\log\frac{\mu_-}{\mu_+}+p+1|^2}
			\bigg)^{\frac 1{p+\frac 12}}
			,
			\\
			& \kappa_3\DE 0,\quad \psi\DE\sqrt{\frac{\kappa_2}{2p+1}}\,\frac 1\pi\log\frac{\mu_-}{\mu_+}
			.
		\end{align*}
	\item If $\mu_+=0$ or $\mu_-=0$, set
		\begin{align*}
			& \kappa_1\DE 2p+1,\quad \kappa_2\DE 0,\quad \kappa_3\DE 0
			,
			\\
			& \psi\DE\Big(\frac 1{2p+1}\Big)^{1-\frac 1{2p+1}}
			\begin{cases}
				-\Big(\mu_+\Gamma(2p+1)^2\Big)^{\frac 1{2p+1}} \CAS \mu_+>0,
				\\[1mm]
				\Big(\mu_-\Gamma(2p+1)^2\Big)^{\frac 1{2p+1}} \CAS \mu_->0
				.
			\end{cases}
		\end{align*}
	\end{Ilist}
	Plugging this data into the formulas of \Cref{M36} leads to the stated assertion.
\end{proof}

\begin{Remark}
\label{M91}
	If one is interested only in the reproducing kernel $K_E$ and not in the function $E$ itself, the formulae from 
	\Cref{M90}\,$(i)$ can be written in a slightly different form. Namely, set 
	\[
		F(z)\DE e^{i\kappa z}M(\alpha+1,2p+1,-2i\kappa z),\quad G(z)\DE e^{i\kappa z}M(\alpha,2p+1,-2i\kappa z)
		,
	\]
	then obviously $A(z)=\frac 12[F(z)+G(z)]$, and \cite[\phantom{}13.4.4]{abramowitz.stegun:1964} shows that 
	\[
		B(z)=\frac{i(2p+1)}{2\kappa}\big[F(z)-G(z)\big]
		.
	\]
	We see that 
	\[
		K_E(z,w)=\frac{i(2p+1)}{2\kappa}\cdot\frac{F(z)G(\ov w)-G(z)F(\ov w)}{z-\ov w}
		.
	\]
\end{Remark}

%---------
%   FINISH
%---------

\printbibliography

\hspace*{-10mm}
\parbox[t]{64mm}{
{\footnotesize
\begin{flushleft}
\rule{64mm}{1pt}\\
	B.\,Eichinger\\
	Institute for Analysis and Scientific Computing\\
	Vienna University of Technology\\
	Wiedner Hauptstra{\ss}e\ 8--10/101\\
	1040 Wien\\
	AUSTRIA\\
	email: \texttt{benjamin.eichinger@asc.tuwien.ac.at}\\[1mm]
	\textit{and}\\[1mm]
	Lancaster University\\
	School of mathematical sciences\\
	LA1 4YF Lancaster\\
	UK\\
	\rule{64mm}{1pt}
\end{flushleft}
}}
\hspace*{6mm}
\parbox[t]{64mm}{
{\footnotesize
\begin{flushleft}
\rule{64mm}{1pt}\\
	H.\,Woracek\\
	Institute for Analysis and Scientific Computing\\
	Vienna University of Technology\\
	Wiedner Hauptstra{\ss}e\ 8--10/101\\
	1040 Wien\\
	AUSTRIA\\
	email: \texttt{harald.woracek@tuwien.ac.at}\\
	\rule{64mm}{1pt}
\end{flushleft}
}}

%---------
%   LABELS
%---------

%%%%%%%%%%%%%%%%%%%%%%%%%%%%%%%%%%%%%%%%%%%%%%%%%%%%%%%%%%%%%%%%%%%%%%%%%%%%%%%%%%%%%%%%%
% disable in final version
%\newpage
%
%\noindent
%{\large\sc labels:}\\[10mm]
%	
%M:\hspace*{3mm}
%\framebox{
%\begin{tabular}{r@{\quad}r@{\quad}r@{\quad}r@{\quad}r@{\quad}r@{\quad}r@{\quad}r@{\quad}r@{\quad}r@{\quad}}
%	** &  1 &  2 &  3 &  4 &  5 &  6 &  7 &  8 &  9 \rule{0pt}{20pt}\kern-5pt\\[2mm]
%	10 & 11 & 12 & 13 & 14 & 15 & 16 & 17 & 18 & 19 \\[2mm]
%	20 & 21 & 22 & 23 & 24 & 25 & 26 & 27 & 28 & 29 \\[2mm]
%	30 & 31 & 32 & 33 & 34 & 35 & 36 & 37 & 38 & 39 \\[2mm]
%	40 & 41 & 42 & 43 & 44 & 45 & 46 & 47 & 48 & 49 \\[2mm]
%	50 & 51 & 52 & 53 & 54 & 55 & 56 & 57 & 58 & 59 \\[2mm]
%	60 & 61 & 62 & 63 & 64 & 65 & 66 & 67 & 68 & 69 \\[2mm]
%	70 & 71 & 72 & 73 & 74 & 75 & 76 & 77 & 78 & 79 \\[2mm]
%	80 & 81 & 82 & 83 & 84 & 85 & 86 & 87 & 88 &  . \\[2mm]
%	90 & 91 & 92 &  . &  . &  . &  . &  . & 98 & 99 \\[2mm]
%	100& 101& 102& 103& 104& 105& 106& 107& 108& 109\\[2mm]
%	110& 111&  . & 113& 114& 115& 116& 117& 118& 119\\[2mm]
%	120&  . &  . & 123& 124& 125& 126& 127& 128&  . \\[2mm]
%	130& 131& 132& 133& 134&  . &  . &  . &  . &  . \\[2mm]
%	 . &  . &  . &  . &  . &  . &  . &  . &  . &  . \\[2mm]
%	 . &  . &  . &  . &  . &  . &  . &  . &  . &  . \\[2mm]
%	 . &  . &  . &  . &  . &  . &  . &  . &  . &  . \\[3mm]
%\end{tabular}
%}
%%%%%%%%%%%%%%%%%%%%%%%%%%%%%%%%%%%%%%%%%%%%%%%%%%%%%%%%%%%%%%%%%%%%%%%%%%%%%%%%%%%%%%%%%

%
%
%

\end{document}